\input amstex
\input amsppt.sty

\magnification=\magstep1 

\vcorrection{-0.8cm}
\def\({\left(}
\def\){\right)}
\def\U{{\tilde U}}
\def\be{\beta}

\catcode`\@=11
\newif\iftab@\tab@false
\newif\ifvtab@\vtab@false
\def\tab{\bgroup\tab@true\vtab@false\vst@bfalse\Strich@false%
   \def\\{\global\hline@@false%
     \ifhline@\global\hline@false\global\hline@@true\fi\cr}
   \edef\l@{\the\leftskip}\ialign\bgroup\hskip\l@##\hfil&&##\hfil\cr}
\def\endtab{\cr\egroup\egroup}
\def\vtab{\vtop\bgroup\vst@bfalse\vtab@true\tab@true\Strich@false%
   \bgroup\def\\{\cr}\ialign\bgroup&##\hfil\cr}
\def\endvtab{\cr\egroup\egroup\egroup}
\def\stab{\D@cke0.5pt\null 
 \bgroup\tab@true\vtab@false\vst@bfalse\Strich@true\Let@@\vspace@
 \normalbaselines\offinterlineskip
  \openup\spreadmlines@
 \edef\l@{\the\leftskip}\ialign
 \bgroup\hskip\l@##\hfil&&##\hfil\crcr}
\def\endstab{\crcr\egroup
 \egroup}
\newif\ifvst@b\vst@bfalse
\def\vstab{\D@cke0.5pt\null
 \vtop\bgroup\tab@true\vtab@false\vst@btrue\Strich@true\bgroup\Let@@\vspace@
 \normalbaselines\offinterlineskip
  \openup\spreadmlines@\bgroup}
\def\endvstab{\crcr\egroup\egroup
 \egroup\tab@false\Strich@false}

\newdimen\htstrut@
\htstrut@8.5\p@
\newdimen\htStrut@
\htStrut@12\p@
\newdimen\dpstrut@
\dpstrut@3.5\p@
\newdimen\dpStrut@
\dpStrut@3.5\p@
\def\openup{\afterassignment\@penup\dimen@=}
\def\@penup{\advance\lineskip\dimen@
  \advance\baselineskip\dimen@
  \advance\lineskiplimit\dimen@
  \divide\dimen@ by2
  \advance\htstrut@\dimen@
  \advance\htStrut@\dimen@
  \advance\dpstrut@\dimen@
  \advance\dpStrut@\dimen@}
\def\Let@@{\relax%
    \def\\{\global\hline@@false%
     \ifhline@\global\hline@false\global\hline@@true\fi\cr}%
    \iffalse}\fi}
\def\matrix{\null\,\vcenter\bgroup
 \tab@false\vtab@false\vst@bfalse\Strich@false\Let@@\vspace@
 \normalbaselines\openup\spreadmlines@\ialign
 \bgroup\hfil$\m@th##$\hfil&&\quad\hfil$\m@th##$\hfil\crcr
 \Mathstrut@\crcr\noalign{\kern-\baselineskip}}
\def\endmatrix{\crcr\Mathstrut@\crcr\noalign{\kern-\baselineskip}\egroup
 \egroup\,}
\def\smatrix{\D@cke0.5pt\null\,
 \vcenter\bgroup\tab@false\vtab@false\vst@bfalse\Strich@true\Let@@\vspace@
 \normalbaselines\offinterlineskip
  \openup\spreadmlines@\ialign
 \bgroup\hfil$\m@th##$\hfil&&\quad\hfil$\m@th##$\hfil\crcr}
\def\endsmatrix{\crcr\egroup
 \egroup\,\Strich@false}
\newdimen\D@cke
\def\Dicke#1{\global\D@cke#1}
\newtoks\tabs@\tabs@{&}
\newif\ifStrich@\Strich@false
\newif\iff@rst

\def\Stricherr@{\iftab@\ifvtab@\errmessage{\noexpand\s not allowed
     here. Use \noexpand\vstab!}%
  \else\errmessage{\noexpand\s not allowed here. Use \noexpand\stab!}%
  \fi\else\errmessage{\noexpand\s not allowed
     here. Use \noexpand\smatrix!}\fi}
\def\format{\ifvst@b\else\crcr\fi\egroup\iffalse{\fi\ifnum`}=0 \fi\format@}
\def\format@#1\\{\def\preamble@{#1}%
 \def\Str@chfehlt##1{\ifx##1\s\Stricherr@\fi\ifx##1\\\let\Next\relax%
   \else\let\Next\Str@chfehlt\fi\Next}%
 \def\c{\hfil\noexpand\ifhline@@\hbox{\vrule height\htStrut@%
   depth\dpstrut@ width\z@}\noexpand\fi%
   \ifStrich@\hbox{\vrule height\htstrut@ depth\dpstrut@ width\z@}%
   \fi\iftab@\else$\m@th\fi\the\hashtoks@\iftab@\else$\fi\hfil}%
 \def\r{\hfil\noexpand\ifhline@@\hbox{\vrule height\htStrut@%
   depth\dpstrut@ width\z@}\noexpand\fi%
   \ifStrich@\hbox{\vrule height\htstrut@ depth\dpstrut@ width\z@}%
   \fi\iftab@\else$\m@th\fi\the\hashtoks@\iftab@\else$\fi}%
 \def\l{\noexpand\ifhline@@\hbox{\vrule height\htStrut@%
   depth\dpstrut@ width\z@}\noexpand\fi%
   \ifStrich@\hbox{\vrule height\htstrut@ depth\dpstrut@ width\z@}%
   \fi\iftab@\else$\m@th\fi\the\hashtoks@\iftab@\else$\fi\hfil}%
 \def\s{\ifStrich@\ \the\tabs@\vrule width\D@cke\the\hashtoks@%
          \fi\the\tabs@\ }%
 \def\sa{\ifStrich@\vrule width\D@cke\the\hashtoks@%
            \the\tabs@\ %
            \fi}%
 \def\se{\ifStrich@\ \the\tabs@\vrule width\D@cke\the\hashtoks@\fi}%
 \def\cd{\hfil\noexpand\ifhline@@\hbox{\vrule height\htStrut@%
   depth\dpstrut@ width\z@}\noexpand\fi%
   \ifStrich@\hbox{\vrule height\htstrut@ depth\dpstrut@ width\z@}%
   \fi$\dsize\m@th\the\hashtoks@$\hfil}%
 \def\rd{\hfil\noexpand\ifhline@@\hbox{\vrule height\htStrut@%
   depth\dpstrut@ width\z@}\noexpand\fi%
   \ifStrich@\hbox{\vrule height\htstrut@ depth\dpstrut@ width\z@}%
   \fi$\dsize\m@th\the\hashtoks@$}%
 \def\ld{\noexpand\ifhline@@\hbox{\vrule height\htStrut@%
   depth\dpstrut@ width\z@}\noexpand\fi%
   \ifStrich@\hbox{\vrule height\htstrut@ depth\dpstrut@ width\z@}%
   \fi$\dsize\m@th\the\hashtoks@$\hfil}%
 \ifStrich@\else\Str@chfehlt#1\\\fi%
 \setbox\z@\hbox{\xdef\Preamble@{\preamble@}}\ifnum`{=0 \fi\iffalse}\fi
 \ialign\bgroup\span\Preamble@\crcr}
\newif\ifhline@\hline@false
\newif\ifhline@@\hline@@false
\def\hlinefor#1{\multispan@{\strip@#1 }\leaders\hrule height\D@cke\hfill%
    \global\hline@true\ignorespaces}
\catcode`\@=13

\catcode`\@=11
\font\tenln    = line10
\font\tenlnw   = linew10

\newskip\Einheit \Einheit=0.5cm
\newcount\xcoord \newcount\ycoord
\newdimen\xdim \newdimen\ydim \newdimen\PfadD@cke \newdimen\Pfadd@cke

\newcount\@tempcnta
\newcount\@tempcntb

\newdimen\@tempdima
\newdimen\@tempdimb

\newdimen\@wholewidth
\newdimen\@halfwidth

\newcount\@xarg
\newcount\@yarg
\newcount\@yyarg
\newbox\@linechar
\newbox\@tempboxa
\newdimen\@linelen
\newdimen\@clnwd
\newdimen\@clnht

\newif\if@negarg

\def\@whilenoop#1{}
\def\@whiledim#1\do #2{\ifdim #1\relax#2\@iwhiledim{#1\relax#2}\fi}
\def\@iwhiledim#1{\ifdim #1\let\@nextwhile=\@iwhiledim
        \else\let\@nextwhile=\@whilenoop\fi\@nextwhile{#1}}

\def\@whileswnoop#1\fi{}
\def\@whilesw#1\fi#2{#1#2\@iwhilesw{#1#2}\fi\fi}
\def\@iwhilesw#1\fi{#1\let\@nextwhile=\@iwhilesw
         \else\let\@nextwhile=\@whileswnoop\fi\@nextwhile{#1}\fi}

\def\thinlines{\let\@linefnt\tenln \let\@circlefnt\tencirc
  \@wholewidth\fontdimen8\tenln \@halfwidth .5\@wholewidth}
\def\thicklines{\let\@linefnt\tenlnw \let\@circlefnt\tencircw
  \@wholewidth\fontdimen8\tenlnw \@halfwidth .5\@wholewidth}
\thinlines

\PfadD@cke1pt \Pfadd@cke0.5pt
\def\PfadDicke#1{\PfadD@cke#1 \divide\PfadD@cke by2 \Pfadd@cke\PfadD@cke \multiply\PfadD@cke by2}
\long\def\LOOP#1\REPEAT{\def\BODY{#1}\ITERATE}
\def\ITERATE{\BODY \let\next\ITERATE \else\let\next\relax\fi \next}
\let\REPEAT=\fi
\def\Punkt{\hbox{\raise-2pt\hbox to0pt{\hss$\ssize\bullet$\hss}}}
\def\DuennPunkt(#1,#2){\unskip
  \raise#2 \Einheit\hbox to0pt{\hskip#1 \Einheit
          \raise-2.5pt\hbox to0pt{\hss$\bullet$\hss}\hss}}
\def\NormalPunkt(#1,#2){\unskip
  \raise#2 \Einheit\hbox to0pt{\hskip#1 \Einheit
          \raise-3pt\hbox to0pt{\hss\twelvepoint$\bullet$\hss}\hss}}
\def\DickPunkt(#1,#2){\unskip
  \raise#2 \Einheit\hbox to0pt{\hskip#1 \Einheit
          \raise-4pt\hbox to0pt{\hss\fourteenpoint$\bullet$\hss}\hss}}
\def\Kreis(#1,#2){\unskip
  \raise#2 \Einheit\hbox to0pt{\hskip#1 \Einheit
          \raise-4pt\hbox to0pt{\hss\fourteenpoint$\circ$\hss}\hss}}

\def\Line@(#1,#2)#3{\@xarg #1\relax \@yarg #2\relax
\@linelen=#3\Einheit
\ifnum\@xarg =0 \@vline
  \else \ifnum\@yarg =0 \@hline \else \@sline\fi
\fi}

\def\@sline{\ifnum\@xarg< 0 \@negargtrue \@xarg -\@xarg \@yyarg -\@yarg
  \else \@negargfalse \@yyarg \@yarg \fi
\ifnum \@yyarg >0 \@tempcnta\@yyarg \else \@tempcnta -\@yyarg \fi
\ifnum\@tempcnta>6 \@badlinearg\@tempcnta0 \fi
\ifnum\@xarg>6 \@badlinearg\@xarg 1 \fi
\setbox\@linechar\hbox{\@linefnt\@getlinechar(\@xarg,\@yyarg)}%
\ifnum \@yarg >0 \let\@upordown\raise \@clnht\z@
   \else\let\@upordown\lower \@clnht \ht\@linechar\fi
\@clnwd=\wd\@linechar
\if@negarg \hskip -\wd\@linechar \def\@tempa{\hskip -2\wd\@linechar}\else
     \let\@tempa\relax \fi
\@whiledim \@clnwd <\@linelen \do
  {\@upordown\@clnht\copy\@linechar
   \@tempa
   \advance\@clnht \ht\@linechar
   \advance\@clnwd \wd\@linechar}%
\advance\@clnht -\ht\@linechar
\advance\@clnwd -\wd\@linechar
\@tempdima\@linelen\advance\@tempdima -\@clnwd
\@tempdimb\@tempdima\advance\@tempdimb -\wd\@linechar
\if@negarg \hskip -\@tempdimb \else \hskip \@tempdimb \fi
\multiply\@tempdima \@m
\@tempcnta \@tempdima \@tempdima \wd\@linechar \divide\@tempcnta \@tempdima
\@tempdima \ht\@linechar \multiply\@tempdima \@tempcnta
\divide\@tempdima \@m
\advance\@clnht \@tempdima
\ifdim \@linelen <\wd\@linechar
   \hskip \wd\@linechar
  \else\@upordown\@clnht\copy\@linechar\fi}

\def\@hline{\ifnum \@xarg <0 \hskip -\@linelen \fi
\vrule height\Pfadd@cke width \@linelen depth\Pfadd@cke
\ifnum \@xarg <0 \hskip -\@linelen \fi}

\def\@getlinechar(#1,#2){\@tempcnta#1\relax\multiply\@tempcnta 8
\advance\@tempcnta -9 \ifnum #2>0 \advance\@tempcnta #2\relax\else
\advance\@tempcnta -#2\relax\advance\@tempcnta 64 \fi
\char\@tempcnta}

\def\Vektor(#1,#2)#3(#4,#5){\unskip\leavevmode
  \xcoord#4\relax \ycoord#5\relax
      \raise\ycoord \Einheit\hbox to0pt{\hskip\xcoord \Einheit
         \Vector@(#1,#2){#3}\hss}}

\def\Vector@(#1,#2)#3{\@xarg #1\relax \@yarg #2\relax
\@tempcnta \ifnum\@xarg<0 -\@xarg\else\@xarg\fi
\ifnum\@tempcnta<5\relax
\@linelen=#3\Einheit
\ifnum\@xarg =0 \@vvector
  \else \ifnum\@yarg =0 \@hvector \else \@svector\fi
\fi
\else\@badlinearg\fi}

\def\@hvector{\@hline\hbox to 0pt{\@linefnt
\ifnum \@xarg <0 \@getlarrow(1,0)\hss\else
    \hss\@getrarrow(1,0)\fi}}

\def\@vvector{\ifnum \@yarg <0 \@downvector \else \@upvector \fi}

\def\@svector{\@sline
\@tempcnta\@yarg \ifnum\@tempcnta <0 \@tempcnta=-\@tempcnta\fi
\ifnum\@tempcnta <5
  \hskip -\wd\@linechar
  \@upordown\@clnht \hbox{\@linefnt  \if@negarg
  \@getlarrow(\@xarg,\@yyarg) \else \@getrarrow(\@xarg,\@yyarg) \fi}%
\else\@badlinearg\fi}

\def\@upline{\hbox to \z@{\hskip -.5\Pfadd@cke \vrule width \Pfadd@cke
   height \@linelen depth \z@\hss}}

\def\@downline{\hbox to \z@{\hskip -.5\Pfadd@cke \vrule width \Pfadd@cke
   height \z@ depth \@linelen \hss}}

\def\@upvector{\@upline\setbox\@tempboxa\hbox{\@linefnt\char'66}\raise
     \@linelen \hbox to\z@{\lower \ht\@tempboxa\box\@tempboxa\hss}}

\def\@downvector{\@downline\lower \@linelen
      \hbox to \z@{\@linefnt\char'77\hss}}

\def\@getlarrow(#1,#2){\ifnum #2 =\z@ \@tempcnta='33\else
\@tempcnta=#1\relax\multiply\@tempcnta \sixt@@n \advance\@tempcnta
-9 \@tempcntb=#2\relax\multiply\@tempcntb \tw@
\ifnum \@tempcntb >0 \advance\@tempcnta \@tempcntb\relax
\else\advance\@tempcnta -\@tempcntb\advance\@tempcnta 64
\fi\fi\char\@tempcnta}

\def\@getrarrow(#1,#2){\@tempcntb=#2\relax
\ifnum\@tempcntb < 0 \@tempcntb=-\@tempcntb\relax\fi
\ifcase \@tempcntb\relax \@tempcnta='55 \or
\ifnum #1<3 \@tempcnta=#1\relax\multiply\@tempcnta
24 \advance\@tempcnta -6 \else \ifnum #1=3 \@tempcnta=49
\else\@tempcnta=58 \fi\fi\or
\ifnum #1<3 \@tempcnta=#1\relax\multiply\@tempcnta
24 \advance\@tempcnta -3 \else \@tempcnta=51\fi\or
\@tempcnta=#1\relax\multiply\@tempcnta
\sixt@@n \advance\@tempcnta -\tw@ \else
\@tempcnta=#1\relax\multiply\@tempcnta
\sixt@@n \advance\@tempcnta 7 \fi\ifnum #2<0 \advance\@tempcnta 64 \fi
\char\@tempcnta}

\def\Diagonale(#1,#2)#3{\unskip\leavevmode
  \xcoord#1\relax \ycoord#2\relax
      \raise\ycoord \Einheit\hbox to0pt{\hskip\xcoord \Einheit
         \Line@(1,1){#3}\hss}}
\def\AntiDiagonale(#1,#2)#3{\unskip\leavevmode
  \xcoord#1\relax \ycoord#2\relax 
      \raise\ycoord \Einheit\hbox to0pt{\hskip\xcoord \Einheit
         \Line@(1,-1){#3}\hss}}
\def\Pfad(#1,#2),#3\endPfad{\unskip\leavevmode
  \xcoord#1 \ycoord#2 \thicklines\ZeichnePfad#3\endPfad\thinlines}
\def\ZeichnePfad#1{\ifx#1\endPfad\let\next\relax
  \else\let\next\ZeichnePfad
    \ifnum#1=1
      \raise\ycoord \Einheit\hbox to0pt{\hskip\xcoord \Einheit
         \vrule height\Pfadd@cke width1 \Einheit depth\Pfadd@cke\hss}%
      \advance\xcoord by 1
    \else\ifnum#1=2
      \raise\ycoord \Einheit\hbox to0pt{\hskip\xcoord \Einheit
        \hbox{\hskip-\PfadD@cke\vrule height1 \Einheit width\PfadD@cke depth0pt}\hss}%
      \advance\ycoord by 1
    \else\ifnum#1=3
      \raise\ycoord \Einheit\hbox to0pt{\hskip\xcoord \Einheit
         \Line@(1,1){1}\hss}
      \advance\xcoord by 1
      \advance\ycoord by 1
    \else\ifnum#1=4
      \raise\ycoord \Einheit\hbox to0pt{\hskip\xcoord \Einheit
         \Line@(1,-1){1}\hss}
      \advance\xcoord by 1
      \advance\ycoord by -1
    \fi\fi\fi\fi
  \fi\next}
\def\hSSchritt{\leavevmode\raise-.4pt\hbox to0pt{\hss.\hss}\hskip.2\Einheit
  \raise-.4pt\hbox to0pt{\hss.\hss}\hskip.2\Einheit
  \raise-.4pt\hbox to0pt{\hss.\hss}\hskip.2\Einheit
  \raise-.4pt\hbox to0pt{\hss.\hss}\hskip.2\Einheit
  \raise-.4pt\hbox to0pt{\hss.\hss}\hskip.2\Einheit}
\def\vSSchritt{\vbox{\baselineskip.2\Einheit\lineskiplimit0pt
\hbox{.}\hbox{.}\hbox{.}\hbox{.}\hbox{.}}}
\def\DSSchritt{\leavevmode\raise-.4pt\hbox to0pt{%
  \hbox to0pt{\hss.\hss}\hskip.2\Einheit
  \raise.2\Einheit\hbox to0pt{\hss.\hss}\hskip.2\Einheit
  \raise.4\Einheit\hbox to0pt{\hss.\hss}\hskip.2\Einheit
  \raise.6\Einheit\hbox to0pt{\hss.\hss}\hskip.2\Einheit
  \raise.8\Einheit\hbox to0pt{\hss.\hss}\hss}}
\def\dSSchritt{\leavevmode\raise-.4pt\hbox to0pt{%
  \hbox to0pt{\hss.\hss}\hskip.2\Einheit
  \raise-.2\Einheit\hbox to0pt{\hss.\hss}\hskip.2\Einheit
  \raise-.4\Einheit\hbox to0pt{\hss.\hss}\hskip.2\Einheit
  \raise-.6\Einheit\hbox to0pt{\hss.\hss}\hskip.2\Einheit
  \raise-.8\Einheit\hbox to0pt{\hss.\hss}\hss}}
\def\SPfad(#1,#2),#3\endSPfad{\unskip\leavevmode
  \xcoord#1 \ycoord#2 \ZeichneSPfad#3\endSPfad}
\def\ZeichneSPfad#1{\ifx#1\endSPfad\let\next\relax
  \else\let\next\ZeichneSPfad
    \ifnum#1=1
      \raise\ycoord \Einheit\hbox to0pt{\hskip\xcoord \Einheit
         \hSSchritt\hss}%
      \advance\xcoord by 1
    \else\ifnum#1=2
      \raise\ycoord \Einheit\hbox to0pt{\hskip\xcoord \Einheit
        \hbox{\hskip-2pt \vSSchritt}\hss}%
      \advance\ycoord by 1
    \else\ifnum#1=3
      \raise\ycoord \Einheit\hbox to0pt{\hskip\xcoord \Einheit
         \DSSchritt\hss}
      \advance\xcoord by 1
      \advance\ycoord by 1
    \else\ifnum#1=4
      \raise\ycoord \Einheit\hbox to0pt{\hskip\xcoord \Einheit
         \dSSchritt\hss}
      \advance\xcoord by 1
      \advance\ycoord by -1
    \fi\fi\fi\fi
  \fi\next}
\def\Koordinatenachsen(#1,#2){\unskip
 \hbox to0pt{\hskip-.5pt\vrule height#2 \Einheit width.5pt depth1 \Einheit}%
 \hbox to0pt{\hskip-1 \Einheit \xcoord#1 \advance\xcoord by1
    \vrule height0.25pt width\xcoord \Einheit depth0.25pt\hss}}
\def\Koordinatenachsen(#1,#2)(#3,#4){\unskip
 \hbox to0pt{\hskip-.5pt \ycoord-#4 \advance\ycoord by1
    \vrule height#2 \Einheit width.5pt depth\ycoord \Einheit}%
 \hbox to0pt{\hskip-1 \Einheit \hskip#3\Einheit 
    \xcoord#1 \advance\xcoord by1 \advance\xcoord by-#3 
    \vrule height0.25pt width\xcoord \Einheit depth0.25pt\hss}}
\def\Gitter(#1,#2){\unskip \xcoord0 \ycoord0 \leavevmode
  \LOOP\ifnum\ycoord<#2
    \loop\ifnum\xcoord<#1
      \raise\ycoord \Einheit\hbox to0pt{\hskip\xcoord \Einheit\Punkt\hss}%
      \advance\xcoord by1
    \repeat
    \xcoord0
    \advance\ycoord by1
  \REPEAT}
\def\Gitter(#1,#2)(#3,#4){\unskip \xcoord#3 \ycoord#4 \leavevmode
  \LOOP\ifnum\ycoord<#2
    \loop\ifnum\xcoord<#1
      \raise\ycoord \Einheit\hbox to0pt{\hskip\xcoord \Einheit\Punkt\hss}%
      \advance\xcoord by1
    \repeat
    \xcoord#3
    \advance\ycoord by1
  \REPEAT}
\def\Label#1#2(#3,#4){\unskip \xdim#3 \Einheit \ydim#4 \Einheit
  \def\lo{\advance\xdim by-.5 \Einheit \advance\ydim by.5 \Einheit}%
  \def\llo{\advance\xdim by-.25cm \advance\ydim by.5 \Einheit}%
  \def\loo{\advance\xdim by-.5 \Einheit \advance\ydim by.25cm}%
  \def\o{\advance\ydim by.25cm}%
  \def\ro{\advance\xdim by.5 \Einheit \advance\ydim by.5 \Einheit}%
  \def\rro{\advance\xdim by.25cm \advance\ydim by.5 \Einheit}%
  \def\roo{\advance\xdim by.5 \Einheit \advance\ydim by.25cm}%
  \def\l{\advance\xdim by-.30cm}%
  \def\r{\advance\xdim by.30cm}%
  \def\lu{\advance\xdim by-.5 \Einheit \advance\ydim by-.6 \Einheit}%
  \def\llu{\advance\xdim by-.25cm \advance\ydim by-.6 \Einheit}%
  \def\luu{\advance\xdim by-.5 \Einheit \advance\ydim by-.30cm}%
  \def\u{\advance\ydim by-.30cm}%
  \def\ru{\advance\xdim by.5 \Einheit \advance\ydim by-.6 \Einheit}%
  \def\rru{\advance\xdim by.25cm \advance\ydim by-.6 \Einheit}%
  \def\ruu{\advance\xdim by.5 \Einheit \advance\ydim by-.30cm}%
  #1\raise\ydim\hbox to0pt{\hskip\xdim
     \vbox to0pt{\vss\hbox to0pt{\hss$#2$\hss}\vss}\hss}%
}
\catcode`\@=13

\catcode`\@=11
\font@\twelverm=cmr10 scaled\magstep1
\font@\twelveit=cmti10 scaled\magstep1
\font@\twelvebf=cmbx10 scaled\magstep1
\font@\twelvei=cmmi10 scaled\magstep1
\font@\twelvesy=cmsy10 scaled\magstep1
\font@\twelveex=cmex10 scaled\magstep1

\newtoks\twelvepoint@
\def\twelvepoint{\normalbaselineskip15\p@
 \abovedisplayskip15\p@ plus3.6\p@ minus10.8\p@
 \belowdisplayskip\abovedisplayskip
 \abovedisplayshortskip\z@ plus3.6\p@
 \belowdisplayshortskip8.4\p@ plus3.6\p@ minus4.8\p@
 \textonlyfont@\rm\twelverm \textonlyfont@\it\twelveit
 \textonlyfont@\sl\twelvesl \textonlyfont@\bf\twelvebf
 \textonlyfont@\smc\twelvesmc \textonlyfont@\tt\twelvett
%
 \ifsyntax@ \def\big##1{{\hbox{$\left##1\right.$}}}%
  \let\Big\big \let\bigg\big \let\Bigg\big
 \else
  \textfont\z@=\twelverm  \scriptfont\z@=\tenrm  \scriptscriptfont\z@=\sevenrm
  \textfont\@ne=\twelvei  \scriptfont\@ne=\teni  \scriptscriptfont\@ne=\seveni
  \textfont\tw@=\twelvesy \scriptfont\tw@=\tensy \scriptscriptfont\tw@=\sevensy
  \textfont\thr@@=\twelveex \scriptfont\thr@@=\tenex
        \scriptscriptfont\thr@@=\tenex
  \textfont\itfam=\twelveit \scriptfont\itfam=\tenit
        \scriptscriptfont\itfam=\tenit
  \textfont\bffam=\twelvebf \scriptfont\bffam=\tenbf
        \scriptscriptfont\bffam=\sevenbf
  \setbox\strutbox\hbox{\vrule height10.2\p@ depth4.2\p@ width\z@}%
  \setbox\strutbox@\hbox{\lower.6\normallineskiplimit\vbox{%
        \kern-\normallineskiplimit\copy\strutbox}}%
 \setbox\z@\vbox{\hbox{$($}\kern\z@}\bigsize@=1.4\ht\z@
 \fi
 \normalbaselines\rm\ex@.2326ex\jot3.6\ex@\the\twelvepoint@}

\font@\fourteenrm=cmr10 scaled\magstep2
\font@\fourteenit=cmti10 scaled\magstep2
\font@\fourteensl=cmsl10 scaled\magstep2
\font@\fourteensmc=cmcsc10 scaled\magstep2
\font@\fourteentt=cmtt10 scaled\magstep2
\font@\fourteenbf=cmbx10 scaled\magstep2
\font@\fourteeni=cmmi10 scaled\magstep2
\font@\fourteensy=cmsy10 scaled\magstep2
\font@\fourteenex=cmex10 scaled\magstep2
\font@\fourteenmsa=msam10 scaled\magstep2
\font@\fourteeneufm=eufm10 scaled\magstep2
\font@\fourteenmsb=msbm10 scaled\magstep2
\newtoks\fourteenpoint@
\def\fourteenpoint{\normalbaselineskip15\p@
 \abovedisplayskip18\p@ plus4.3\p@ minus12.9\p@
 \belowdisplayskip\abovedisplayskip
 \abovedisplayshortskip\z@ plus4.3\p@
 \belowdisplayshortskip10.1\p@ plus4.3\p@ minus5.8\p@
 \textonlyfont@\rm\fourteenrm \textonlyfont@\it\fourteenit
 \textonlyfont@\sl\fourteensl \textonlyfont@\bf\fourteenbf
 \textonlyfont@\smc\fourteensmc \textonlyfont@\tt\fourteentt
%
 \ifsyntax@ \def\big##1{{\hbox{$\left##1\right.$}}}%
  \let\Big\big \let\bigg\big \let\Bigg\big
 \else
  \textfont\z@=\fourteenrm  \scriptfont\z@=\twelverm  \scriptscriptfont\z@=\tenrm
  \textfont\@ne=\fourteeni  \scriptfont\@ne=\twelvei  \scriptscriptfont\@ne=\teni
  \textfont\tw@=\fourteensy \scriptfont\tw@=\twelvesy \scriptscriptfont\tw@=\tensy
  \textfont\thr@@=\fourteenex \scriptfont\thr@@=\twelveex
        \scriptscriptfont\thr@@=\twelveex
  \textfont\itfam=\fourteenit \scriptfont\itfam=\twelveit
        \scriptscriptfont\itfam=\twelveit
  \textfont\bffam=\fourteenbf \scriptfont\bffam=\twelvebf
        \scriptscriptfont\bffam=\tenbf
  \setbox\strutbox\hbox{\vrule height12.2\p@ depth5\p@ width\z@}%
  \setbox\strutbox@\hbox{\lower.72\normallineskiplimit\vbox{%
        \kern-\normallineskiplimit\copy\strutbox}}%
 \setbox\z@\vbox{\hbox{$($}\kern\z@}\bigsize@=1.7\ht\z@
 \fi
 \normalbaselines\rm\ex@.2326ex\jot4.3\ex@\the\fourteenpoint@}

\font@\seventeenrm=cmr10 scaled\magstep3
\font@\seventeenit=cmti10 scaled\magstep3
\font@\seventeensl=cmsl10 scaled\magstep3
\font@\seventeensmc=cmcsc10 scaled\magstep3
\font@\seventeentt=cmtt10 scaled\magstep3
\font@\seventeenbf=cmbx10 scaled\magstep3
\font@\seventeeni=cmmi10 scaled\magstep3
\font@\seventeensy=cmsy10 scaled\magstep3
\font@\seventeenex=cmex10 scaled\magstep3
\font@\seventeenmsa=msam10 scaled\magstep3
\font@\seventeeneufm=eufm10 scaled\magstep3
\font@\seventeenmsb=msbm10 scaled\magstep3
\newtoks\seventeenpoint@
\def\seventeenpoint{\normalbaselineskip18\p@
 \abovedisplayskip21.6\p@ plus5.2\p@ minus15.4\p@
 \belowdisplayskip\abovedisplayskip
 \abovedisplayshortskip\z@ plus5.2\p@
 \belowdisplayshortskip12.1\p@ plus5.2\p@ minus7\p@
 \textonlyfont@\rm\seventeenrm \textonlyfont@\it\seventeenit
 \textonlyfont@\sl\seventeensl \textonlyfont@\bf\seventeenbf
 \textonlyfont@\smc\seventeensmc \textonlyfont@\tt\seventeentt
%
 \ifsyntax@ \def\big##1{{\hbox{$\left##1\right.$}}}%
  \let\Big\big \let\bigg\big \let\Bigg\big
 \else
  \textfont\z@=\seventeenrm  \scriptfont\z@=\fourteenrm  \scriptscriptfont\z@=\twelverm
  \textfont\@ne=\seventeeni  \scriptfont\@ne=\fourteeni  \scriptscriptfont\@ne=\twelvei
  \textfont\tw@=\seventeensy \scriptfont\tw@=\fourteensy \scriptscriptfont\tw@=\twelvesy
  \textfont\thr@@=\seventeenex \scriptfont\thr@@=\fourteenex
        \scriptscriptfont\thr@@=\fourteenex
  \textfont\itfam=\seventeenit \scriptfont\itfam=\fourteenit
        \scriptscriptfont\itfam=\fourteenit
  \textfont\bffam=\seventeenbf \scriptfont\bffam=\fourteenbf
        \scriptscriptfont\bffam=\twelvebf
  \setbox\strutbox\hbox{\vrule height14.6\p@ depth6\p@ width\z@}%
  \setbox\strutbox@\hbox{\lower.86\normallineskiplimit\vbox{%
        \kern-\normallineskiplimit\copy\strutbox}}%
 \setbox\z@\vbox{\hbox{$($}\kern\z@}\bigsize@=2\ht\z@
 \fi
 \normalbaselines\rm\ex@.2326ex\jot5.2\ex@\the\seventeenpoint@}

\catcode`\@=13

\topmatter
\author  
Christian Krattenthaler\footnote"$^\dagger$"{\hbox{Partially supported by the Austrian
Science Foundation FWF, grant P13190-MAT.}}\\
Luigi Orsina\\ Paolo Papi\endauthor
\bigskip
\address
\vskip 0.pt Institut f\"ur Mathematik der Universit\"at Wien,
\vskip 0.pt Strudlhofgasse 4, A-1090 Wien, Austria.
\vskip 0.pt e-mail: {KRATT\@Ap.Univie.Ac.At}
\vskip 3pt
\vskip 5.pt Dipartimento di Matematica, Istituto G. Castelnuovo
\vskip 0.pt Universit\`a di Roma ``La Sapienza"  
\vskip 0.pt Piazzale Aldo Moro 5
\vskip 0.pt 00185 Rome --- ITALY 
\vskip 0.pt e-mail:{\rm \ orsina\@mat.uniroma1.it,\hskip3 pt papi\@mat.uniroma1.it}
 \endaddress

\leftheadtext {C. Krattenthaler, L. Orsina and P. Papi}
\rightheadtext { Enumeration of $ad$-nilpotent
$\frak b$-ideals for simple Lie algebras}
\abstract
We  provide explicit formulas for the number of $ad$-nilpotent ideals of a Borel subalgebra of a complex simple Lie
algebra  having fixed class of nilpotence.
\endabstract
\title  Enumeration of $ad$-nilpotent
$\frak b$-ideals for simple Lie algebras
\endtitle
\keywords ad-nilpotent ideal, Lie algebra, order ideal, Dyck path,
Chebyshev polynomial\endkeywords
\subjclass Primary: 17B20; Secondary: 05A15 05A19 05E15 17B30\endsubjclass
\endtopmatter
\nologo
\def\endemo{\qed\enddemo}

\def\g{\frak g}
\def\h{\frak h}
\def\n{\frak n}
\def\bb{\frak b}
\def\D{\Delta}
\def\l{\lambda}
\def\Dp{\Delta^+}

\def\d{\delta}
\def\r{R(\h)}

\def\a{\alpha}
\def\b{\beta}

\def\I{\Cal I^n}
\def\l{\lambda}
\def\d{\delta}

\def\i{{\frak i}}

\def\p{\Phi}

\TagsOnRight

\def\l{\lambda}
\def\v#1{\left\vert#1\right\vert}
\def\fl#1{\left\lfloor#1\right\rfloor}
\def\cl#1{\left\lceil#1\right\rceil}

\def\({\left(}
\def\){\right)}

\document

\heading\S1 Introduction\endheading

In this paper we  provide  formulas for the number of $ad$-nilpotent ideals of a Borel subalgebra of a complex simple Lie
algebra
$\g$ having 
fixed class of nilpotence. Our results continue and 
complete the analysis of
\cite{OP} and \cite{AKOP}, where Lie algebras of
type $A$ have been treated.

The basic framework is the following.
Let $\g$ be a  complex simple Lie algebra of rank $n$. Let
$\h\subset\g$ be a fixed Cartan subalgebra,  
$\D$ the corresponding  root system of $\g$. Fix a positive system
$\Dp$ in $\D$.
For each $\a\in \Dp$ let ${\frak g}_\a$ be the root space of $\g$ relative to $\a$, 
$\n = \bigoplus \limits_{\a\in \Dp}{\frak g}_\a$,  and $\bb$ be the Borel subalgebra
$\bb=\h\oplus \n$.

Let $\I$ denote the set of $ad$-nilpotent ideals (i.e., consisting of
$ad$-nilpotent elements) of $\bb$. Then $\i\in \I$ if and only if $\i$ decomposes as $\i=\bigoplus \limits_{\a\in \p}{\frak
g}_\a$, with $\p\subseteq \Dp$ being a dual order ideal 
in the poset $(\Dp,\leq)$ of positive roots (i.e.,  a subset of the
positive roots with the
property that whenever $\a\ge \be$ with $\be\in\Dp$, then also
$\a\in\Dp$; the partial order $\le$ being the restriction to $\Dp$ of the usual
partial order on the root lattice). Clearly, an ideal
$\i\in
\I$ is nilpotent. We denote its 
class of nilpotence by $n(\i)$. 
(By definition, the class of nilpotence of an
ideal $\i$ is the smallest number
$m$ such that $m$-fold bracketing of $\i$ with itself gives the zero ideal. 
Thus, Abelian ideals are exactly those with class of nilpotence
at most 1.)

$ad$-Nilpotent ideals were studied by Kostant in \cite{K1}, \cite{K2} in connection with representation theory 
of compact semisimple Lie groups. In particular, in \cite{K2}, using ideas of Peterson, Kostant found an encoding of the
Abelian  ideals by means of elements of the affine
Weyl group of $\g$. Moreover, he proved that the 
Abelian ideals are $2^n$ in number (an observation originally made by
Peterson).
These results lead to ask similar questions for $ad$-nilpotent ideals: find an appropriate encoding and 
enumerate $\I$ by class of nilpotence. The first problem was settled in \cite{CP}, where an encoding in the spirit of 
Kostant--Peterson was found. The second question, for $\g$ of type
$A$, was addressed and given a solution in \cite{OP} and \cite{AKOP}.
Interestingly, it turns out that in type $A$ ideals in $\I$ with class 
of nilpotence $K$ are equinumerous with Dyck paths of length $2n+2$
and of 
height $K+1$ (see \cite{AKOP, Theorem~4.4}). In fact, in
\cite{AKOP, Sec.~5} an explicit bijection between the former set of
ideals and the latter set of Dyck paths is constructed.
Several explicit formulas for the number of these ideals can be given 
(see \cite{OP}
and \cite{AKOP, Theorems~4.2, 4.5, 4.6}).

In the present paper we resolve the enumerative question in all other
cases. In types $B$, $C$, and $D$ we give explicit formulas, expressed 
in terms of Chebyshev polynomials of the second kind, for the
generating function for ideals in $\I$ with class of nilpotence $K$
(see the theorems labelled B, C, and D in Section~2).
In types $E$, $F$, and $G$ we have computed the number of ideals in
$\I$ with a given class of nilpotence in a rather straightforward way, 
using a computer (see Section~8).

In type $C$ we are again able to establish a relationship between the
ideals in $\I$ and certain paths: in Theorem~C.3 we prove that in type
$C$ ideals in $\I$ with class of nilpotence $K$ are equinumerous with 
paths with step vectors $(1,1)$ and $(1,-1)$ of length $2n$ and height 
$K+1$, which start at the origin and never pass
below the $x$-axis. However, in contrast to type $A$, we are not
able to provide an explicit bijection between this set of ideals and
this set of paths. 
In types $B$ and $D$ we do not even know any relationship between the
ideals in $\I$  
with class of nilpotence $K$ and any set of paths. It seems indeed that paths 
are not the right objects to consider in the case of
type $B$ and type $D$. We believe that a fruitful direction for further
investigation is to try to relate ideals in $\I$ in any of the
types $A$, $B$, $C$, or $D$ to non-crossing partitions of the
corresponding types. Indeed, in any of these types, the sequence of numbers 
of ideals in $\I$ is identical with the sequence of numbers of 
non-crossing partitions (cf\. \cite{Si} and
\cite{R, Cor.~10}). We speculate that
there is a naturally defined statistics on non-crossing partitions,
such that, in any of the types $A$, $B$, $C$, or $D$, the ideals in
$\I$ with class of nilpotence $K$ are equinumerous with the
non-crossing partitions (in the set corresponding to $\I$)
on which this statistics has the value $K$. (In fact, a rougher
version of that problem was already stated in Remark~2 of \cite{R}.
There, the problem is posed to relate non-crossing partitions to {\it
non-nesting} partitions. However, non-nesting partitions are by definition
antichains of roots, which in turn are in bijection
with our $ad$-nilpotent ideals. See also the final paragraph of the
introduction.)

Our paper is organized as follows. In the next section we first
summarize the enumerative results from \cite{OP} and \cite{AKOP} in
type $A$, in order to make a comparison with our results in types $B$, 
$C$, and $D$ easy. We state the results in types $B$, $C$ and $D$
subsequently, also in Section~2. 
In Section~3 we review 
the encoding of $ad$-nilpotent $\frak b$-ideals in type $A$ in terms of
subdiagrams of the staircase shape, and the geometric algorithm to
determine class of nilpotence in type $A$ given in
\cite{AKOP}. This algorithm is of particular importance for our
analysis in types $B$, $C$, and $D$. For in Section~4 we show that
$ad$-nilpotent $\frak b$-ideals in types $B$, $C$, and $D$ can be encoded as
certain subdiagrams of a {\it shifted\/} staircase
(as done originally by Shi \cite{Sh}), and, for the determination of class of  
nilpotence, we may in fact apply the type $A$ algorithm to diagrams
that are obtained by conjugation operations from the subdiagrams of
the shifted staircase. We use these observations to prove the
enumerative results of Section~2 in the subsequent sections. Section~5 
is devoted to the proofs in type $C$, Section~6 to the proofs in type
$B$, and Section~7 to the proofs in type $D$. While the results follow
quite easily in type $C$, types $B$ and $D$ need a more refined (though
still elementary) analysis. Finally, our enumerative results in the
exceptional types are reported in Section~8.

It is interesting to note that $ad$-nilpotent ideals in a Borel
subalgebra have appeared in various (equivalent) disguises in the
literature. In \cite{Sh} Shi studied {\it $\oplus$-sign types}, 
a relevant notion in the theory of {\it cells} for
Weyl groups. As he shows
these are in bijection with dual order ideals in the set of positive
roots, which in turn are in bijection with our $ad$-nilpotent ideals,
as we have argued right at the beginning. In another direction, Postnikov
introduced, for any of the classical and exceptional root systems,
{\it non-nesting partitions} (see \cite{R, Remark~2}). 
As we already mentioned, these are by definition
antichains of positive roots (a definition that
translates to a nice combinatorial
interpretation in the classical types), which in turn are in bijection
with dual order ideals in the set of positive roots, and, hence, 
again with our $ad$-nilpotent ideals. Postnikov also showed that
non-nesting partitions are in bijection with the regions of the
{\it Catalan arrangement\/} contained in the fundamental chamber. 
These hyperplane arrangements have for instance been
investigated by Athanasiadis \cite{A1, A2}.
Interestingly, Reiner \cite{R,
Remark~2} states a uniform formula for the
total number of non-nesting partitions, and, thus, for the {\it
total number} of $ad$-nilpotent $\frak b$-ideals in a Borel
subalgebra for $\g$ of any type, namely
$$\prod_{i=1}^n\frac {h+e_i+1} {e_i+1},\tag1.1$$
where $h$ is the Coxeter number  
and $e_1,e_2,\ldots, e_n$ are the exponents of $\D$ (cf\. \cite{Hu;
Sec.~3.19}). This formula can easily be checked case by
case, by comparing it with the known results on these numbers.
For the classical types one may consult \cite{CP,
Theorem~3.1} or \cite{Sh, Sec.~3}, and for the exceptional
types \cite{Sh, Sec.~3}, except that the number for $E_6$ in \cite{Sh} 
has to be corrected; see the
last paragraph and Table~1 in Section~8. Alternatively, for the
classical types, one may
consult \cite{A1, Theorem~5.5} (cf\. also \cite{A2, Theorem~3.1}).
However, it would be desirable to find a uniform proof for this
formula. 

\head \S2 The enumeration of $ad$-nilpotent $\frak b$-ideals in the
classical types\endhead

We begin by recalling the enumerative results from \cite{OP} and
\cite{AKOP} in type $A$.

\proclaim{Theorem A.1}\cite{OP}, \cite{AKOP, Theorem~4.2}  
Let $\g$ be of type $A_n$.  Let $\a_n(K)$ denote the number of ideals in
$\I $ with class of nilpotence $K$. Then 
$$
\a_n(K)=\sum_{0=i_0<i_1<\dots<i_K<i_{K+1}=n+1}\prod_{j=0}^{K-1}
\binom{i_{j+2}-i_j-1}{i_{j+1}-i_j}.\tag A.1$$\endproclaim

The combinatorial meaning of the previous formula has been  analysed
in detail in \cite{AKOP}. In particular, it connects to the
enumeration of Dyck paths. Recall that a Dyck path
is a lattice path with diagonal step vectors
$(1,1)$ and $(1,-1)$, which starts at the origin and returns to the
$x$-axis, and which does not pass below the $x$-axis. 
We define the height
of a Dyck path to be the maximum ordinate of its peaks.

\proclaim{Theorem A.2}\cite{AKOP, Theorem~4.4}
For $\g$ of type $A_n$, 
the number of ideals in $\I$ with class of nilpotence $K$ 
is exactly the same as the number of
Dyck paths from $(0,0)$ to $(2n+2,0)$ with height $K+1$.
\endproclaim

Generating function results for Dyck paths translate into the
following result for ideals with a given class of nilpotence (cf\. the
fact quoted in the paragraph before Theorem~C.3). Let
$U_n(x)$ denote the $n$th Chebyshev 
polynomial of the second 
kind, $U_n(\cos t)=\sin((n+1)t)/\sin t$, or, explicitly,
$$U_n(x)=\sum_{j\geq 0}(-1)^j\binom{n-j}{j}(2x)^{n-2j}.$$
In the statement of the following result, and also later in the
paper, we write $\U_k$ as short-hand for
$U_k\(1/2\sqrt x\)$. 

\proclaim{Theorem A.3}\cite{AKOP, Theorem~4.6} 
Let $\g$ be of type $A_n$. 
Let $\a^\le_n(h)$ denote the number of ideals in $\I$ with class of
nilpotence at most $h$. Then
$$\align
1+\sum_{n=0}^\infty\a^\le_n(h) x^{n+1}&=
\frac{\U_{h+1}}{\sqrt x\,\U_{h+2}}.
\tag A.3 \endalign$$\endproclaim

A nice implication of this theorem is that not only is the number of
Abelian ideals given by a nice compact formula, namely $2^n$, but
there are also similarly nice formulae for ideals with class of
nilpotence at most $2$
and $3$.

\proclaim{Corollary A.4} \cite{AKOP, Corollary~4.7}
For $\g$ of type $A_n$, 
the number of ideals in $\I$ with class of nilpotence
at most $2$ is the Fibonacci number $F_{2n}$.
The number of ideals in $\I$ with class of nilpotence
at most $3$ is $(3^n+1)/2$. 
\endproclaim

Denote by 
 $\a_{n}(h,K)$ the number of such ideals with
dimension $h$ and class
of nilpotence $K$ and set
$$C_n(q,t)=\sum_{h,K\geq 0}\a_n(h,K)t^hq^K.\tag A.4$$
This notation reflects the fact that $C_n(1,1)$ is the $(n+1)$st Catalan
number\linebreak $\frac{1}{n+2}\binom{2n+2}{n+1}$, so that the polynomials $C_n(q,t)$ are
$(q,t)$-analogues of the Catalan numbers\footnote{Although these
$(q,t)$-analogues of the Catalan numbers are different from the
$(q,t)$-Catalan numbers of Garsia and Haiman \cite{GHn}, a
combinatorial interpretation of the latter $(q,t)$-Catalan numbers
has recently been found by Haglund \cite{H} (the proof \cite{GHd}
jointly with Garsia) which is amazingly similar to the combinatorial
interpretation of (A.4) of \cite{AKOP} (recalled in Section~2) in
terms of two statistics on Ferrers diagrams contained in a staircase.
Whereas the first statistics is (basically) identical
in both cases (in \cite{AKOP} it is the area of the diagram,
respectively in     
\cite{H} and \cite{GHd} it is the difference of the area of the
staircase and the 
area of the diagram), the difference in the second statistics can be best
described by saying that the statistics of \cite{AKOP} (the number of
touching points on the diagonal line that bounds the Ferrers diagram;
see Figure~1) is a descent-like statistics, whereas Haglund's
statistics is the corresponding major-like statistics.}. 
In \cite{AKOP} we proved the following formula for $C_n(q,t)$. 
\proclaim{Theorem A.5}\cite{AKOP, Theorem 6.1} We have
$$C_n(q,t)=\kern-1pt\sum_{K=0}^n
\left(\sum_{0=i_0<i_1<\dots<i_K<i_{K+1}=n+1} 
\prod_{j=0}^{K-1}t^{i_{j+1}(i_{j+3}-i_{j+2})}
\bmatrix {i_{j+2}-i_j-1}\\{i_{j+1}-i_j}\endbmatrix_t
\right)q^K,\tag A.5$$
with $i_{K+2}=n+2$. Here,
$\left[\smallmatrix{m}\\{n}\endsmallmatrix\right]_t$ is 
the $t$-binomial coefficient, defined by
$$\bmatrix{m}\\{n}\endbmatrix_t=\cases
\frac{(1-t^m)(1-t^{m-1})\cdots(1-t^{m-n+1})}
{(1-t^n)(1-t^{n-1})\cdots(1-t)} \quad & \text{if $m \geq n > 0$,} \\
1 \quad & \text{if $n=0$,}\\
0 \quad & \text{in any other case.} \endcases
$$
\endproclaim

\medskip
In this paper we provide analogues of Theorems~A.1 and A.3 and of
Corollary~A.4 
for $\g$ of type $B$, $C$, and $D$. We are able to give an analogue of
Theorem~A.2 only in type $C$. Finally, 
again in type $C$, we also provide an analogue of
Theorem~A.5, thus   
obtaining a $(q,t)$-analogue of the binomial coefficient
$\binom{2n}{n}$. In principle, we would also be able to write down
analogues of Theorem~A.5 in types $B$ and $D$. However, the
corresponding formulas are rather unwieldy, 
so that we prefer to omit these
for the sake of brevity.

\medskip
Now we state these 
results. We begin with the type $C$ analogue of Theorem~A.1.

\proclaim{Theorem C.1} Let $\g$ be of type $C_n$. Let $\gamma_n(K)$
denote the number of ideals in $\I$ with class of nilpotence
$K$. Then 
$$\gamma_n(K)=\cases
\sum_{0<i_1<\dots<i_k<i_{k+1}=n} 
\prod_{j=1}^{k-1}
\binom{i_{j+2}-i_j-1}{i_{j+1}-i_j}
\cdot
\sum _{\ell=0} ^{i_2-i_1-1}\binom{i_1+i_2-1}{\ell}\\
&\hskip-1.3cm\text{if $K=2k$,}\\
\sum_{-i_2<i_1\le 0<i_2<\dots<i_k<i_{k+1}=n} 
2^{i_1+i_2-1}
\prod_{j=1}^{k-1}
\binom{i_{j+2}-i_j-1}{i_{j+1}-i_j}\\
&\hskip-1.3cm\text{if $K=2k-1$,}\endcases\tag C.1$$
where $i_{K+2}=n+1$. 
\endproclaim

Next we state the type $C$ analogue of Theorem~A.3.

\proclaim{Theorem C.2} 
Let $\g$ be of type $C_n$. 
Let $\gamma^\le_n(h)$ be the number of ideals in
$\I$ with class of nilpotence at 
most $h$. Then the generating function 
$\sum _{n \ge0} ^{}\gamma^\le_n(h)x^n$ is given by
$$\dfrac {\sum _{i=0} ^{\cl{h/2}}\U_{h+1-2i}} {\sqrt x\,\U_{h+2}}=
\dfrac {\U_{h+1}+\U_{h-1}+\U_{h-3}+\cdots} {\sqrt x\,\U_{h+2}}\ .
\tag C.2$$
\endproclaim

Theorem~C.2 allows to relate the ideals of given class of nilpotence 
again to path combinatorics. Appealing to
the fact (see for example \cite{Kr, Theorem~A2 + Fact~A3})
that the generating function $\sum _{P}
^{}x^{\ell(P)/2}$ for paths $P$ with step vectors $(1,1)$ and $(1,-1)$, 
which start at the origin, never
exceed height $k$ (the height being again defined as the maximum
ordinate of the peaks of the path), and end
at height $s$ ($\ell(P)$ denoting the length of $P$) is given by
$$\frac {\U_{k-s}} {\sqrt x\,\U_{k+1}},$$
we obtain the following type $C$ analogue of
Theorem~A.2. 

\proclaim{Theorem C.3}
For $\g$ of type $C_n$, 
the number of ideals in $\I$ with class of nilpotence $K$ 
is exactly the same as the number of
paths with step vectors $(1,1)$ and $(1,-1)$, which start at the origin, 
have length $2n$ and height $K+1$, and
never pass below the $x$-axis. 
\endproclaim

As we already said in the introduction, it would be desirable to find
an explicit bijection between the ideals and the paths in the above theorem.

This lattice path interpretation of the number of ideals with given
class of nilpotence, together with 
the iterated reflection principle formula for paths
that stay between two parallel lines (cf\. \cite{Mo, Ch.~1, Th.~2}), 
implies,  upon little
simplification, an explicit formula for this number.
This formula must be preferred over the one in (C.1),
as it is much simpler and computationally superior. (An analogous
formula in type $A$ has been given in \cite{AKOP, Eq.~(4.6)}.)

\proclaim{Corollary~C.4}
For $\g$ of type $C_n$, 
the number of ideals in $\I$ with class of nilpotence at most $h$ 
is given by
$$\sum _{s=0} ^{\fl{h/2}}\sum _{k\in\Bbb Z} ^{}\frac {1+2s+2k(h+2)}
{2n+1} \binom {2n+1}{n-s-k(h+2)}.\tag C.4$$
\endproclaim

Using Theorem~C.2 again, 
it is easy to deduce the following type $C$ analogue of Corollary~A.4. 
There, $F_n$ denotes again the $n$th Fibonacci number.
\proclaim{Corollary C.5}
Let $\g$ be of type $C_n$. Then for $n\ge1$
the number of ideals in $\I$ with class of
nilpotence at most $2$ is $F_{2n}$, and 
the number of ideals with class of nilpotence at most $3$ is $2\cdot 3^{n-1}$.
\endproclaim

Finally, we state the type $C$ analogue of Theorem~A.5. 
\proclaim{Theorem C.6} 
Let $\g$ be of type $C_n$.
Let $\gamma_n(h,K)$ be the number of ideals in $\I$ with
dimension $h$ and class of nilpotence $K$. Then
$$\multline
\sum _{h,K\ge0} ^{}\gamma_n(h,K)t^hq^{K}=\sum _{k=0} ^{n}
\sum_{-i_2<i_1<i_2<\dots<i_k<i_{k+1}=n} q^{2k-\chi(i_1\le0)}\\
\cdot \bigg(
\prod_{j=1}^{k-1}t^{(i_{j+1}+n)(i_{j+3}-i_{j+2})}
\bmatrix {i_{j+2}-i_j-1}\\{i_{j+1}-i_j}\endbmatrix_t\\
\cdot
\sum _{\ell=0} ^{i_2-i_1-1}\bmatrix i_1+i_2-1\\\ell\endbmatrix_t
t^{(i_1+n)(i_3-i_2)-\binom{n-i_2+1}2+\binom {\ell+1}2}\bigg),
\endmultline\tag C.6$$
where $i_{k+2}=n+1$, 
$\left[\smallmatrix{m}\\{n}\endsmallmatrix\right]_t$ is the
$t$-binomial coefficient defined in the statement of Theorem~A.5,
and $\chi$ is the usual truth function, $\chi(\Cal A)=1$ if $\Cal A$ is
true and $\chi(\Cal A)=0$ otherwise.\endproclaim
Clearly, since in type $C$ the total number of ideals in $\I$ is the
central binomial coefficient $\binom{2n}{n}$ (see Corollary~4.2), 
the expression in (C.6) is a $(q,t)$-analogue of $\binom{2n}{n}$.

\remark{Remark}The innermost sum in (C.6) does in fact simplify if
$i_1\le0$. For, by means of the $q$-binomial theorem (cf\. \cite{GR,
Ex.~1.2({\it vi})}), we have
$$\align \sum _{\ell=0} ^{i_2-i_1-1}\bmatrix i_1+i_2-1\\\ell\endbmatrix_t
t^{\binom {\ell+1}2}
&=\sum _{\ell=0} ^{\infty}\bmatrix i_1+i_2-1\\\ell\endbmatrix_t
t^{\binom {\ell+1}2}\\
&=(1+t)(1+t^2)\cdots(1+t^{i_1+i_2-1}). 
\tag C.7\endalign$$
In view of this, it is now obvious that the special case $t=1$ of Theorem~C.6
implies Theorem~C.1.
\endremark

\medskip
Now we state our results for type $B$. 

\proclaim{Theorem B.1}
Let $\g$ be of type $B_n$. 
Let $\be_n(K)$ denote the number of 
ideals in $\I$ with class of nilpotence $K$.
Then the generating function $\sum _{n \ge0} ^{}\be_n(K)x^n$ is given by
$$\cases \dsize \frac {\U_{2k}+\U_k\U_{k+1}\U_{2k-1}} {\sqrt
x\,\U_{2k} \U_{2k+1}\U_{2k+2}}&\text{if }K=2k,\\
\dsize \frac {\U_{2k}+\U_{k+1}^2\U_{2k-2}+\U_{k-1}^2\U_{2k-2}+\U_2} {
\U_{2k-1} \U_{2k}\U_{2k+1}}&\text{if }K=2k-1.
\endcases\tag B.1$$
\endproclaim
Analogues of formula (A.1) for type $B$ will be established in the
proof of the above theorem in Section~6. 
An easy consequence of Theorem~B.1 is the following type $B$ analogue
of Theorem~A.3. It can be proved by a straightforward induction on $h$.
\proclaim{Corollary B.2}
Let $\g$ be of type $B_n$. 
Let $\be^\le_n(h)$ denote the number of ideals in $\I$ with class of
nilpotence at most $h$. Then the
generating function $\sum _{n \ge0} ^{}\be^\le_n(h)x^n$ is given by 
$$\cases
  {\dfrac{\sum_{i = 1}^{{h /2}}
        \left( 2\,i + 1 \right) \,\U_{2i - 1} + 
      ( h +1) \U_{h+1} + \sum_{i = 1}^{{h /2}}
        2i\,\U_{2h + 3 - 2i}  
      } {\U_{h+1}\,\U_{h+2}}}&\text{if $h$ is even,}\\
  {\dfrac{\sum_{i = 1}^{{(h - 1)/2}}
        \left( 2\,i + 1 \right) \,\U_{2i - 1} + 
      ( h +1) \U_h + \sum_{i = 1}^{{(h + 1)/2}}
        2i\,\U_{2h + 3 - 2i} 
      } {\U_{h+1}\,\U_{h+2}}}&\text{if $h$ is odd.}
\endcases\tag B.2$$\endproclaim

The above theorem implies a type $B$
analogue of Corollary~A.4.

\proclaim{Corollary B.3}
Let $\g$ be of type $B_n$. 
Then for $n\ge1$
the number of ideals in $\I$ with class of nilpotence 
at most $2$ is $F_{2n}+F_{2n-2}-2^{n-1}$, and 
the number of ideals with class of nilpotence at most $3$ is 
$\frac {1} {2}(5\cdot 3^{n-1}+1)-F_{2n-2}$.
\endproclaim

Next we state our results for type $D$.

\proclaim{Theorem D.1}
Let $\g$ be of type $D_n$. 
Let $\d_n(K)$ denote the number of ideals in $\I$ with
class of nilpotence $K$. 
Then the generating function $\sum _{n \ge0} ^{}\d_n(K)x^n$ is given by
$x/(1-x)$ if $K=0$, and otherwise by
$$\cases \dsize \frac {2-\U_{2k}+2\U_{2k+2}+3\U_{k}\U_{k+1}\U_{2k-1}} {\sqrt
x\,\U_{2k} \U_{2k+1}\U_{2k+2}}
&\hskip-2cm\text{if }K=2k,\\
\dsize \frac {2\U_{2k+2}+\U_{2k}-3\U_{2k-2}+\U_{k}\U_{k+1}\U_{2k-1}+4\U_{k}^2\U_{2k-2}+
             \U_{k-1}\U_{k}\U_{2k-3}+2} {
\U_{2k-1} \U_{2k}\U_{2k+1}}\\
&\hskip-2cm\text{if }K=2k-1.
\endcases\tag D.1$$
\endproclaim
Analogues of formula (A.1) for type $D$ will be established in the
proof of the above theorem in Section~7. 
An easy consequence of Theorem~D.1 is the following type $D$ analogue of
Theorem~A.3. Once again, it can be proved by a straightforward
induction on $h$. 
\proclaim{Corollary D.2}
Let $\g$ be of type $D_n$. 
Let $\d^\le_n(h)$ denote the number of ideals in $\I$ 
with class of nilpotence at most $h$. Then the
generating function $\sum _{n \ge0} ^{}\d^\le_n(h)x^n$ is given by
$x/(1-x)$ if $h=0$, by $x(1+2x)/(1-2x)$ if $h=1$, and otherwise by
$$\cases
  x\Big(6\U_{1}+\sum_{i = 1}^{{(h-2) /2}}
        \left( 6i + 8 \right) \,\U_{2i + 1} + 
      ( 3h +4) \U_{h+1} \\
   \quad + \sum_{i = 0}^{{(h-2) /2}}
        (6i+5)\,\U_{2h + 1 - 2i}  +\U_{2h+3}\Big)\Big/
       {\U_{h+1}\U_{h+2}}\quad &\text{if $h$ is even,}\\
  x\Big(6\U_{1}+\sum_{i = 1}^{{(h - 3)/2}}
        \left( 6i + 8 \right) \,\U_{2i + 1} + 
      ( 3h +4) \U_h \\
    \quad + \sum_{i = 0}^{{(h - 1)/2}}
        (6i+5)\,\U_{2h + 1 - 2i} +\U_{2h+3}\Big)\Big/
       {\U_{h+1}\U_{h+2}}\quad &\text{if $h$ is odd.}
\endcases\tag D.2$$\endproclaim

The above theorem implies a type $D$
analogue of Corollary~A.4.

\proclaim{Corollary D.3}
Let $\g$ be of type $D_n$. 
Then for $n\ge2$
the number of ideals in $\I$ with class of nilpotence 
at most $2$ is $5F_{2n-3}-2^{n-2}$, and 
the number of ideals with class of nilpotence at most $3$ is 
$\frac {13} {2}\cdot 3^{n-2}-\frac {3} {2}+4F_{2n}-7F_{2n-1}$.
\endproclaim

\heading\S3 $ad$-nilpotent ideals in type $A$ \endheading

Throughout this section, $\g$ will be of type $A_n$, i.e., 
$\g$ is the Lie
algebra $sl(n+1,\Bbb C)$ of $(n+1)\times (n+1)$ traceless matrices.
In the following paragraphs we collect the findings from \cite{AKOP} 
on how to
encode ideals in $\I$ and how to efficiently compute their class of 
nilpotence.

Recall from the introduction that an ideal $\i$ is in $\I$, i.e., is an
$ad$-nilpotent ideal of our fixed Borel subalgebra $\bb$, if and only
if $\i$ can be written as $\i=\bigoplus \limits_{\a\in \p}{\frak
g}_\a$, where $\p\subseteq \Dp$ is a dual order ideal 
in the poset $(\Dp,\leq)$ of positive roots.
This allows us to represent
$ad$-nilpotent ideals conveniently in a geometric fashion, which will
be crucial in all subsequent considerations. Clearly, any positive
root in $A_n$ can be written as a sum of simple roots. Explicitly,  
with $\a_1$, $\a_2$, \dots, $\a_n$ denoting the simple roots, let us write
$\tau_{i,j}=\a_i+\dots+\a_{n-j+1},$  $1\leq i\le n$, $1\le j\leq
n-i+1$. If we place the roots $\tau_{i,j}$, $j=1,2,\dots,n-i+1$, in the
$i$th row of a diagram, then this defines an arrangement of
the positive roots in a staircase fashion. For example, for
$A_3$ we obtain the arrangement
$$\matrix\format\l\quad &\l\quad &\l\\
\a_1+\a_2+\a_3 &\a_1+\a_2&\a_1\\
\a_2+\a_3 &\a_2 & \\
\a_3 \endmatrix$$
Obviously, the above defines an identification of the positive
roots with the cells of the staircase diagram $(n,n-1,\dots,1)$ (for
all partition notation we refer the reader to \cite{Ma, Ch.~I, Sec.~1}), 
in which the root $\tau_{i,j}$
is identified with the cell $(i,j)$. 
Given an $ad$-nilpotent ideal $\i$, written as 
$\i=\bigoplus\limits_{\a\in \p}{\frak g}_\a$, for some collection $\p$ of
positive roots, we can use the above
identification to represent $\i$ as the set of cells that corresponds
to the roots in $\p$. Since, as we noted above, $\i$ is a dual order
ideal, the set of cells obtained forms a (Ferrers diagram of a)
partition. For example, the ideal ${\frak g}_{\a_1+\a_2+\a_3}\oplus
{\frak g}_{\a_1+\a_2}\oplus {\frak g}_{\a_1}\oplus {\frak g}_{\a_2+\a_3}$ corresponds to the
partition $(3,1,0)$. 
This defines a bijection between $ad$-nilpotent ideals in $sl(n+1,\Bbb
C)$ and subdiagrams of $(n,n-1,\dots,1)$.

Now, let $m_{i,j}$ be the maximal
number $m$ such that the root space ${\frak g}_{\tau_{i,j}}$ occurs in 
$$\i^{m}:=\underbrace {[\cdots[\i,\i],\dots]}_{m\text{ occurrences of
}\i}.$$
Alternatively, $m_{i,j}$ is the maximal number $m$ such that there is
a decomposition of $\tau_{i,j}$ of the form $\tau_{i,j}=\b_1+\dots+\b_m$
with all $\b_i$'s in $\p$, the collection of positive roots corresponding to
$\i$. Written succinctly, 
$$m_{i,j}=\max \{m:\tau_{i,j}=\b_1+\dots+\b_m,\text{ for some
}\b_l\in\p,\ l=1,2,\dots,m\}.\tag3.1$$

By definition, the class of nilpotence of $\i$  is given by 
$\max_{i,j} m_{i,j}$. 
We claim, first, that the numbers $m_{i,j}$ satisfy $m_{i,j}\ge m_{i+1,j}$
and $m_{i,j}\ge m_{i,j+1}$ for all $i$ and $j$, and, second, 
that the $m_{i,j}$ (hence, in particular, the class of
nilpotence) can be obtained by the following algorithm (see \cite{AKOP,
Prop.~3.1}):
let $\l$ be the subdiagram of $(n,n-1,\dots,1)$ that corresponds to
$\i$ according to the identification explained above. Define a
filling $(t_{i,j})_{1\le i\le n,\,1\le j\le n-i+1}$ 
of the cells of $(n,n-1,\dots,1)$ by recursively setting
$$t_{i,j}=\cases 0\quad &\text{if\ } (i,j)\notin \lambda,\\1\quad &\text{if\
} (i,j)\text{ is a corner cell of }\lambda,\\
\max\limits_{j<k\leq n-i+1}\ \{t_{i,k}+t_{n-k+2,j}\}\quad
&\text{otherwise.}\endcases
\tag 3.2 $$ 
It is easy to see that the above rule uniquely defines a filling of
the shape $(n,n-1,\dots,1)$, whose nonzero entries are precisely those
corresponding to the cells of $\l$. E.g., when $n=4$, the fillings 
corresponding to 
$(2,1,0,0)$, $(3,3,2,1)$, $(4,3,2,1)$ are respectively
$$
\matrix
1&1&0&0\\1&0&0\\0&0\\0\endmatrix\qquad
\matrix
3&2&1&0\\3&2&1\\2&1\\1\endmatrix\qquad
\matrix
4&3&2&1\\3&2&1\\2&1\\1\endmatrix$$

The claim is that $m_{i,j}=t_{i,j}$ for all $i$ and $j$.

Let us go through a proof of these two facts that allows to be ``recycled"
when we shall discuss the computation of class of nilpotence in types
$B$, $C$, and $D$. For proving the first claim, suppose that there is
a representation of $\tau_{i+1,j}$ as 
$$\tau_{i+1,j}=\b_1+\dots+\b_m,\tag3.3$$
with all $\b_l$'s in $\p$. Then one of the $\b_l$'s must be equal to
some $\tau_{i+1,k}$, because otherwise the sum on the right-hand side
of (3.3), when expanded as a sum of simple roots, would either
contain no $\a_{i+1}$ or in addition to $\a_{i+1}$ also $\a_{i}$.
Without loss of generality let $\b_1=\tau_{i+1,k}$. Then
$$\tau_{i,j}=\a_i+\tau_{i+1,j}=\a_i+\b_1+\dots+\b_m=
\tau_{i,k}+\b_2+\dots+\b_m.$$
Hence we have $m_{i,j}\ge m_{i+1,j}$. A similar argument proves
$m_{i,j}\ge m_{i,j+1}$.

For proving the second claim,
we do a reverse induction on $i+j$. The claim is obvious if $i+j=n+1$.
For the induction step, assume that the claim is right for all
$(i',j')$ with $i'+j'>i+j$. Now consider the root $\tau_{i,j}$. If
$j<k\le n-i+1$, then we can write it as
$$\tau_{i,j}=\tau_{i,k}+\tau_{n-k+2,j}.\tag3.4$$ 
By the induction hypothesis
and (3.1) we know that 
$$t_{i,k}=\max \{m:\tau_{i,k}=\b_1+\dots+\b_m,\text{ for some
}\b_l\in\p,\ l=1,2,\dots,m\}$$
and
$$t_{n-k+2,j}=\max \{m:\tau_{n-k+2,j}=\b_1+\dots+\b_m,\text{ for some
}\b_l\in\p,\ l=1,2,\dots,m\}\tag3.5$$
Combined with (3.4), this implies that 
$$t_{i,k}+t_{n-k+2,j}\le
m_{i,j}\tag3.6$$ 
as long as $t_{i,k}$ and $t_{n-k+2,j}$ are nonzero. But (3.6) is also
true if one or both of $t_{i,k}$ and $t_{n-k+2,j}$ should be zero. If
both are zero, then there is nothing to prove. If, for instance,
$t_{i,k}$ is nonzero, then we have $t_{i,k}=m_{i,k}\le m_{i,j}$, the
equality being true by induction hypothesis, the inequality being
true because of our first claim.

But $m_{i,j}$ cannot be larger than $t_{i,k}+t_{n-k+2,j}$. 
For, as we already noted (just replace $i$ by $i-1$ in the previous
argument), in any decomposition
$$\tau_{i,j}=\b_1+\dots+\b_m,\tag3.7$$
with the $\b_l$'s being positive roots, one of the $\b_l$'s must be
equal to some $\tau_{i,k}$. Again, without loss of generality let
$\b_1=\tau_{i,k}$. Then we have
$\b_2+\dots+\b_m=\tau_{i,j}-\tau_{i,k}=\tau_{n-k+2,j}$. By induction
hypothesis we have (3.5), and hence $m-1\le t_{n-k+2,j}$. It follows
that $m\le 1+t_{n-k+2,j}\le t_{i,k}+t_{n-k+2,j}\le t_{i,j}$, as required.

Since $m_{i,j}\ge m_{i+1,j}$ and $m_{i,j}\ge m_{i,j+1}$, the same
property must be true for the $t_{i,j}$'s.
In particular, this implies that the maximum of all entries is entry
$t_{1,1}$, so that in fact the class of nilpotence of an
ideal $\i$ is equal to $t_{1,1}$.

\medskip
If one is only interested to quickly compute the class of nilpotence
of some ideal $\i$ (i.e., just the entry $t_{1,1}$ of the filling
constructed by (3.2)), then there is a short-cut
through the algorithm (3.2). For
a convenient statement of the result, we write, in abuse of notation,
$n(\l_1,\l_2,\ldots,\l_n)$ for $n(\i)$, given that the 
partition corresponding to $\i$ is $(\l_1,\l_2,\ldots,\l_n)$.

\proclaim{Proposition 2.1}\cite{AKOP, Prop.~3.2}
Let $\i\in\I$ and let 
$\l=(\l_1,\ldots,\l_n)$ be the corresponding partition.
We have $n(0,0,\dots,0)=0$, and otherwise
$$n(\l_1,\l_2,\ldots,\l_n)=n(\l_{n+2-\l_1},\ldots,\l_n)+1.\tag 3.8$$
\endproclaim

It should be noted that on the left-hand side of (3.8) there appears
the class of nilpotence of an ideal in $\I$, whereas on the right-hand
side there appears the class of nilpotence of an ideal in $\Cal
I^{\l_1-1}$ (with corresponding partition
$(\l_{n+2-\l_1},\ldots,\l_n)$). 
The computation, however, can be carried out completely formally,
without reference to ideals, as we now demonstrate by an example.

\example{Example}
Let $\i\in\Cal I^{13}$ be the ideal which corresponds to the partition 
$(10,10,9,6,\mathbreak 5,4,4,3,1,1,1,1,0)$. (This is the partition in
Figure~1. At this point, all dotted lines should be ignored.) 
Then, by applying Proposition~2.1
iteratively, we obtain for the class of nilpotence of $\i$:
$$\align
n(\i)&=n(10,10,9,6,5,4,4,3,1,1,1,1,0)\\ 
&= n(5,4,4,3,1,1,1,1,0)+1\\
&= n(1,1,1,0)+2\\
&=3.
\endalign$$
\endexample

As is obvious from the example, iterated application of
Proposition~2.1 provides a very efficient algorithm for the
determination of the
class of nilpotence of a given ideal $\i$. 

Since it will be essential subsequently, 
we wish to point out that this algorithm has a very nice geometric
rendering. Let, as before, $\l=(\l_1,\l_2,\dots,\l_n)$ 
be the partition corresponding to
$\i$. Consider the Ferrers diagram of $\l$. As it is contained in the
staircase diagram $(n,n-1,\dots,1)$, it must not cross the
antidiagonal line $x+y=n+1$.
We draw a zig-zag line
as follows (see Figure~1, where $n=13$ and
$\l=(10,10,9,6,5,4,4,3,1,1,1,1,0)$): 
we start on the vertical edge on the right of cell
$(1,\l_1)$, and move downward until we touch the antidiagonal
$x+y=n+1$. At the touching point we turn direction from 
vertical-down to horizontal-left, and move on until we touch a vertical
part of the Ferrers diagram. At the touching point we turn direction
from  horizontal-left to vertical-down. Now the procedure is
iterated, until we reach the line $x=0$. The class of nilpotence of
the ideal $\i$ is equal to the number of touching points on
$x+y=n+1$. In Figure~1, the resulting zig-zag line is the dotted line
outside the Ferrers diagram of $(10,10,9,6,5,4,4,3,1,1,1,1,0)$.
There are three touching points on $x+y=n+1=14$, in accordance with
$n(10,10,9,6,5,4,4,3,1,1,1,1,0)=3$, as we computed just above.

\midinsert
\vskip10pt
\vbox{
$$
\Pfad(0,1),2222222222221111111111\endPfad
\Pfad(0,1),1222211212212121112122\endPfad
\SPfad(0,0),1222\endSPfad
\SPfad(1,4),11112222\endSPfad
\SPfad(5,9),11111222\endSPfad
\thinlines
\Diagonale(0,-1){14}
\Label\r{x+y=n+1}(15,12)
\DuennPunkt(1,0)
\DuennPunkt(5,4)
\DuennPunkt(10,9)
\hskip7cm
$$
\centerline{\eightpoint Figure 1}
}
\vskip10pt
\endinsert

\heading\S4 $ad$-nilpotent ideals in types $B, C, D$ \endheading

Let $\g$ be a  complex simple Lie algebra of  type  $B_n,\,C_n,$ or
$D_n$. In this section we describe a  diagrammatic encoding of the 
positive roots (cf\. \cite{Sh, Sec.~2--3}) similar 
to the one introduced in the previous section for type $A$.
Throughout this section, for any of the three types, the vectors
$\a_1,\a_2,\dots,\a_n$ denote a basis of simple roots corresponding
to the chosen system $\Dp$ of positive roots.

We arrange the positive roots in a shifted staircase diagram of shape
$(2n-1,2n-3,\ldots,1)$ for $B_n$ and
$C_n$, and of shape $(2n-2,2n-4,\ldots,2)$ for $D_n$ 
as follows.

If $\g$ is of type $C_n$, then we associate the cell $(i,j)$, $1\le
i\le j\le 2n-i$, with the positive root
$$\cases \a_i+\dots+\a_{j-1}+2(\a_{j}+\dots+\a_{n-1})+\a_n\qquad & \text{if }
j\leq n-1,\\ 
\a_i+\dots+\a_{2n-j}\qquad &\text{if }n\leq j\leq 2n-i.\endcases\tag4.C$$
For example, if $n=3$, this defines the following arrangement of
positive roots,
$$\alignat5
&2\a_1+2\a_2+\a_3\qquad &&\a_1+2\a_2+\a_3\qquad &&\a_1+\a_2+\a_3\qquad &&\a_1+\a_2\qquad &&\a_1\\
& &&2\a_2+\a_3 &&\a_2+\a_3 &&\a_2 && \\
& && &&\a_3 &&  &&\endalignat$$

Likewise, if $\g$ is of type $B_n$, then we associate the cell $(i,j)$, $1\le
i\le j\le 2n-i$, with the positive root
$$\cases \a_i+\dots+\a_{j}+2(\a_{j+1}+\dots+\a_n)\qquad &  \text{if }j\leq n-1,\\
\a_i+\dots+\a_{2n-j}\qquad &\text{if }n\leq j\leq 2n-i.\endcases\tag4.B$$
For example, if $n=3$, this defines the following arrangement of
positive roots,
$$\alignat5
&\a_1+2\a_2+2\a_3\qquad &&\a_1+\a_2+2\a_3\qquad &&\a_1+\a_2+\a_3\qquad &&\a_1+\a_2\qquad &&\a_1\\
& &&\a_2+2\a_3 &&\a_2+\a_3 &&\a_2 && \\
& && &&\a_3 &&  &&\endalignat$$

If $\g$ is of type $D_n$, then we associate the cell $(i,j)$, $1\le
i\le j\le 2n-1-i$, with the positive root
$$\cases \a_i+\dots+\a_j+2(\a_{j+1}+\dots+\a_{n-2})+\a_{n-1}+\a_n& \text{if
} j\leq n-2,\\
\a_i+\dots+\a_{n-2}+\a_n&\text{if }j=n-1,\\
\a_i+\dots+\a_{2n-j-1} &\text{if }n\leq j\leq 2n-1-i.\endcases\tag4.D$$
For example, if $n=4$, this defines the following arrangement of
positive roots,
$$\eightpoint\alignat6
&\a_1+2\a_2+\a_3+\a_4\kern15pt &&\a_1+\a_2+\a_3+\a_4\kern15pt&&\a_1+\a_2+\a_4\kern15pt &&\a_1+\a_2+\a_3\kern15pt &&\a_1+\a_2\kern15pt&&\a_1\\
 & &&\a_2+\a_3+\a_4 &&\a_2+\a_4 &&\a_2+\a_3&&\a_2 && \\
& && &&\a_4&&\a_3 &&  &&\endalignat$$

Now, for $\g$ of any of the types $B_n$, $C_n$, or $D_n$, we 
again associate an ideal $\i\in\I$, represented as before as 
$\i=\bigoplus\limits_{\a\in \p}{\frak g}_\a$, to a collection of cells, via the
above identification of positive roots with cells. The resulting
collections are characterized in the proposition below. For
convenience, if $A\subseteq (2n-2,2n-4,\ldots,2)$ is a collection of cells,
we denote by $A^\bullet$ the collection of cells obtained from $A$ by
switching columns $n-1$ and $n$. By a subdiagram of some {\it shifted\/} 
staircase we mean as usual a collection of cells contained in the
shifted staircase which forms a {\it shifted\/} Ferrers diagram
(cf\. \cite{Ma, Ch.~I, Sec.~1, Ex.~9}). 
\proclaim{Proposition 4.1}
If $\g$ is of type $B_n$ or $C_n$, then, under the above identification,
the ideals in $\I$ correspond bijectively to the
subdiagrams of $(2n-1,2n-3,\dots,1)$.
If $\g$ is of type $D_n$, the above identification defines a 
bijection between the ideals in $\I$ 
and the collection of cells $A\subseteq (2n-2,2n-4,\dots,2)$ 
such that either $A$ or $A^\bullet$ is a subdiagram of
$(2n-2,2n-4,\dots,2)$. 
\endproclaim

\proclaim{Corollary 4.2}
The number of ideals in $\I$ is given by
$$
|\I|=\cases\binom{2n}{n}\qquad &\text{ if $\g$ is of type $B_n$ or $C_n$,}\\
\binom{2n}{n}-\binom{2n-2}{n-1}\qquad &\text{ if $\g$ is of type $D_n$.}\endcases$$
\endproclaim

\remark{Remark} 
The previous proposition and corollary are also proved in 
\cite{CP, Theorem~3.1} and, in equivalent form, in \cite{Sh, Sec.~3}. 
In particular,
Corollary~4.2 follows almost immediately from Proposition~4.1. The
most direct argument (which is different from the ones in \cite{CP}
and \cite{Sh}) is as follows:
for $\g$ of type $B_n$ or $C_n$, one maps the shifted
subdiagrams of $(2n-1,2n-3,\dots,1)$ to lattice paths which start at the
origin, never pass below the $x$-axis, and consist of
$2n$ steps, the steps being
up-steps $(1,1)$ and down-steps $(1,-1)$, as before. 
The correspondence between
subdiagrams and paths is best explained with an example at hand. Let
$n=9$ and consider the shifted partition $(16,13,11,8,7,5,3)$, see
Figure~4. Rotate the figure by $45^\circ$ in the positive direction
and then flip it across a vertical line. Then the zig-zag line which
forms the (right) border of the shifted partition becomes such a
lattice path (on attaching a few up-steps at the beginning of the
path). See Figure~2 for the path corresponding to the
partition of Figure~4. (To make a comparison easy, the steps which
correspond to thick segments in Figure~4 are also made thick in Figure~2.)

\midinsert
\vskip10pt
\vbox{
$$
\Gitter(19,5)(0,0)
\Koordinatenachsen(19,5)(0,0)
\Pfad(0,0),334334343344343433\endPfad
\hbox{\hskip1.2pt}
\Pfad(2,2),4\endPfad
\hbox{\hskip-2.4pt}
\Pfad(2,2),4\endPfad
\hbox{\hskip2.4pt}
\Pfad(7,3),4\endPfad
\hbox{\hskip-2.4pt}
\Pfad(7,3),4\endPfad
\hbox{\hskip2.4pt}
\Pfad(13,3),4\endPfad
\hbox{\hskip-2.4pt}
\Pfad(13,3),4\endPfad
\hskip9cm
$$
\centerline{\eightpoint Figure 2}
}
\vskip10pt
\endinsert

Now one can resort to the well-known result (see e.g\. \cite{F,
Theorem~1, (4.6) in Ch.~III, Sec.4}) that the number of these lattice
paths is $\binom {2n}n$. The result for type $D_n$ is then an easy
consequence. (A very different proof can be found in \cite{CP}, where
an enumeration result for trapezoidal plane partitions due to Proctor
is used.)
It can be checked that in all cases the formulas are in accordance with
formula~(1.1).
\endremark

\medskip
We explain now in a rough form 
how to calculate the class of nilpotence by means of
the algorithm used in Section~3 for type $A$. For convenience, let us write
$S_{2n-1}$ for the shifted staircase $(2n-1,2n-3,\dots,1)$, 
$S_{2n-2}$ for the shifted staircase $(2n-2,2n-4,\dots,2)$, and
$T_N$ for the {\it ordinary} staircase diagram $(N,N-1,\dots,1)$. The general 
procedure will be as follows: we  embed the shifted staircase
$S_{N}$, where $N=2n-1$ or $N=2n-2$, respectively, 
into the ordinary staircase $T_{N}$, in such a way
that to any collection of cells $A\subseteq S_N$ corresponds uniquely
a collection of cells $\tilde A\subseteq T_{N}$. 
We shall eventually prove in Proposition~4.3 that the class of nilpotence of an ideal
$\i\in\I$ can be obtained by applying the algorithm of Section~3 to
$\tilde A$, respectively to $\widetilde {A^\bullet}$ in the case of
type $D_n$ if
$A$ is not a diagram, where $A$ is the collection of cells
correponding to $\i$ in the way that was explained earlier.

To explain this in detail, let first $\g$ be of type $C_n$. 
We fill the ordinary staircase $T_{2n-1}$ with the positive
roots of $C_n$ by associating the cells of the ``upper
half" of $T_{2n-1}$ with positive roots as in (4.C), 
and by associating the cells in the ``lower half" of
$T_{2n-1}$ in a symmetric fashion in such a way that cell $(i,j)$ gets
associated with the same root as cell $(j,i)$. (I.e., the line
formed by the cells $(i,i)$, $i=1,2,\dots,n$, constitutes a symmetry
axis of this arrangement of the positive roots.)
Thus, if $n=3$, this defines the
arrangement
$$\alignat5
&2\a_1+2\a_2+\a_3\qquad &&\a_1+2\a_2+\a_3\qquad &&\a_1+\a_2+\a_3\qquad &&\a_1+\a_2\qquad &&\a_1\\
&\a_1+2\a_2+\a_3 &&2\a_2+\a_3 &&\a_2+\a_3 &&\a_2 && \\
&\a_1+\a_2+\a_3 &&\a_2+\a_3 &&\a_3\\
&\a_1+\a_2 &&\a_2\\
&\a_1\endalignat$$
Now, given an ideal $\i\in\I$, written as 
$\i=\bigoplus\limits_{\a\in \p}{\frak g}_\a$, 
with associated collection of cells
$A\subseteq S_{2n-1}$ (which is in fact of the form of a 
shifted Ferrers diagram), 
we let $\tilde A$ be the collection of all cells in $T_{2n-1}$ that
contain a root
$\a\in\p$. Phrased differently, $\tilde A$ is the union of $A$ with its
mirror image about the symmetry axis. Hence, $\tilde A$ forms a
self-conjugate (cf\. \cite{Ma, Ch.~I, Sec.~1}) partition. For example, according
to this description, the ideal 
$\g_{2\a_1+2\a_2+\a_3}+\g_{\a_1+2\a_2+\a_3}+\g_{\a_1+\a_2+\a_3}+
\g_{2\a_2+\a_3}$
corresponds to the self-conjugate partition $(3,2,1)$.

If $\g$ is of type $B_n$, 
we fill the ordinary staircase $T_{2n-1}$ with elements 
of the root lattice of $B_n$ by associating the cells of the ``upper
half" of $T_{2n-1}$ with positive roots as in (4.B), by associating the cell
$(i,i-1)$ with the element of the root lattice
$2\a_i+2\a_{i+1}+\dots+2\a_{n}$ (this element is
not a root!), $i=2,3,\dots,n$, 
and by associating the cells in the ``lower half" of
$T_{2n-1}$ in a symmetric fashion in such a way that cell $(i,j)$ gets
associated with the same root as cell $(j+1,i-1)$. (I.e., here the line
formed by the cells $(i,i-1)$, $i=2,3,\dots,n$, constitutes a symmetry
axis of this arrangement of elements of the root lattice.)
Thus, if $n=3$, this defines the
arrangement
$$\alignat5
&\a_1+2\a_2+2\a_3\qquad &&\a_1+\a_2+2\a_3\qquad &&\a_1+\a_2+\a_3\qquad &&\a_1+\a_2\qquad &&\a_1\\
&2\a_2+2\a_3 &&\a_2+2\a_3 &&\a_2+\a_3 &&\a_2 && \\
&\a_2+2\a_3 &&2\a_3 &&\a_3\\
&\a_2+\a_3 &&\a_3\\
&\a_2\endalignat$$
Here, given an ideal $\i\in\I$, written as 
$\i=\bigoplus\limits_{\a\in \p}{\frak g}_\a$, 
with associated collection of cells
$A\subseteq S_{2n-1}$ (which is in fact of the form of a 
shifted Ferrers diagram), 
we let $\tilde A$ be the collection of all cells in $T_{2n-1}$ that
contain a root $\a\in\p$, 
together with all cells $(i,i-1)$ that are in the same row as
some cell of $A$. Phrased differently, $\tilde A$ is the union of $A$ with its
mirror image about the line formed by the cells $(i,i-1)$,
$i=2,3,\dots,n$, including the cells on that line which are necessary
to ``fill the holes". Hence, if we write $\tilde A$ as a
partition $\l=(\l_1,\l_2,\dots,\l_{2n-1})$, then $(\l_2,\dots,\l_{2n-1})$ is a
self-conjugate partition. Moreover, we have $\l_i\ne i-1$ for $i\ge2$.
For example, according
to this description, the ideal 
$\g_{\a_1+2\a_2+2\a_3}+\g_{\a_1+\a_2+2\a_3}+\g_{\a_2+2\a_3}$
corresponds to the partition $(2,2,1)$.

Finally, if $\g$ is of type $D_n$, 
we fill the ordinary staircase $T_{2n-2}$ with elements 
of the root lattice of $D_n$ in a similar fashion as in type $B$. To
be precise, we associate the cells of the ``upper
half" of $T_{2n-2}$ with positive roots as in (4.D), we associate cell
$(i,i-1)$ with the element of the root lattice
$2\a_i+2\a_{i+1}+\dots+2\a_{n-2}+\a_{n-1}+\a_{n}$
(this element is 
not a root!), $i=2,3,\dots,n-1$, 
and by associating the cells in the ``lower half" of
$T_{2n-1}$ in a symmetric fashion in such a way that cell $(i,j)$ gets
associated with the same root as cell $(j+1,i-1)$. (I.e., the line
formed by the cells $(i,i-1)$, $i=2,3,\dots,n$, constitutes again a symmetry
axis of this arrangement of elements of the root lattice.)
In particular, we do not associate the cell $(n,n-1)$ with any element
of the root lattice.
Thus, if $n=4$, this defines the
arrangement (there, a cross indicates the cell that is not associated
to anything)
$$\eightpoint\alignat6
&\a_1+2\a_2+\a_3+\a_4\kern15pt &&\a_1+\a_2+\a_3+\a_4\kern15pt&&\a_1+\a_2+\a_4\kern15pt &&\a_1+\a_2+\a_3\kern15pt &&\a_1+\a_2\kern15pt&&\a_1\\
 &2\a_2+\a_3+\a_4 &&\a_2+\a_3+\a_4 &&\a_2+\a_4 &&\a_2+\a_3&&\a_2 && \\
&\a_2+\a_3+\a_4 &&\a_3+\a_4 &&\a_4&&\a_3\\
&\a_2+\a_4 &&\a_4&&\times\\
&\a_2+\a_3 &&\a_3\\
&\a_2\endalignat$$
Given an ideal $\i\in\I$, with associated collection of cells
$A\subseteq S_{2n-2}$, we form $\tilde A$ in the same way as in type
$B_n$, i.e., $\tilde A$ is the the union of $A$ with its
mirror image about the line formed by the cells $(i,i-1)$,
$i=2,3,\dots,n-1$, including the cells on that line which are in the
same row as some cell in $A$. Should $A$ not be a shifted diagram 
(and, hence, $\tilde A$ not be an ordinary diagram), then we shall
interchange the $(n-1)$st and the $n$th column and the $n$th and
$(n+1)$st row, and call the resulting (ordinary) diagram
$\widetilde{A^\bullet}$.
Whichever of $\tilde A$ or $\widetilde{A^\bullet}$ is a diagram, 
it is contained in $T_{2n-2}$, and
when written as a partition
$(\l_1,\l_2,\dots,\l_{2n-2})$, it has the property that
$(\l_2,\dots,\l_{2n-2})$ is a
self-conjugate partition and that $\l_i\ne i-1$ for $i\ge2$.
For example, according
to this description, the set $\widetilde{A^\bullet}$ corresponding to
the ideal  
$$\g_{\a_1+2\a_2+\a_3+\a_4}+\g_{\a_1+\a_2+\a_3+\a_4}+\g_{\a_1+\a_2+\a_4}+
\g_{\a_1+\a_2+\a_3}+\g_{\a_2+\a_3+\a_4}+\g_{\a_2+\a_3}$$
is the partition $(4,3,1,1)$.

After ``having set the scene," we come to the crucial result of this
section. As promised earlier, it allows to compute the class of
nilpotence of ideals in any of the types $B_n$, $C_n$, or $D_n$, by applying
the algorithm of Section~3 to a suitable partition.
\proclaim{Proposition 4.3} Let $\g$ be of type $B_n$, $C_n$, or
$D_n$. Let $\i$ be in $\I$, and let $\l$ be the corresponding
{ordinary} partition according to the case-by-case description given
above, i.e., either it stands for $\tilde A$ or for
$\widetilde{A^\bullet}$, and, in addition, $\l=(\l_1,\l_2,\dots)$ 
is self-conjugate in
the case of type $C_n$, whereas in types $B_n$ and $D_n$ it is
$(\l_2,\dots)$ which is self-conjugate. 
Then $n(\i)=n(\l)$, where $n(\l)$ stands for
the result of the algorithm of Proposition~2.1.
\endproclaim 
\remark{Remark}
We should clarify what is meant by ``the result of the algorithm of
Proposition~2.1." For type $B_n$ and type $C_n$ the 
partition $\l$ is contained in $T_{2n-1}$. Accordingly we write
$\l=(\l_1,\dots,\l_{2n-1})$ (i.e., with $2n-1$ components) and apply
(3.8) with $n$ replaced by $2n-1$, and then iterate. For type $D_n$,
however, the 
partition $\l$ is contained in $T_{2n-2}$. Thus in this case we write
$\l=(\l_1,\dots,\l_{2n-2})$ (i.e., with $2n-2$ components) and apply
(3.8) with $n$ replaced by $2n-2$, and then iterate. In the geometric
rendering of the algorithm (see Figure~1), the touching points of the
broken ray are on the line $x+y=2n$ for type $B_n$ and $C_n$, and on
$x+y=2n-1$ for type $D_n$. 
\endremark
\demo{Proof} Construct the filling $(t_{i,j})$ corresponding to $\l$
according to the algorithm (3.2). Let $\tau_{i,j}$ denote the element
of the root lattice associated with cell $(i,j)$. 

Now we copy the arguments of Section~3 that show that $t_{i,j}$
is equal to the right-hand side of (3.1). Since (3.4) is also valid
for our $\tau_{i,j}$'s (regardless whether we consider $B_n$, $C_n$,
or $D_n$!), everything runs through smoothly, with just two modifications.
First, the collection of roots $\p$ (which define the ideal $\i$) must
be replaced by the elements of the root lattice contained in the cells of
$A$. (The reader must recall that in types $B_n$ and $D_n$ this includes
some nonroots.)
Second, if $i>j$, the argument that in any decomposition (3.3) of
$\tau_{i+1,j}$, 
with the $\b_l$'s being positive roots, one of the $\b_l$'s must be
equal to some $\tau_{i+1,k}$, must now be modified by maintaining that
one of the $\b_l$'s must be equal to some $\tau_{h,j}$  
(because, otherwise, the sum on the right-hand side
of (3.3) would either contain no $\a_{j+1}$ or in addition to $\a_{j+1}$ also
$\a_{j}$). Clearly, if $i=j$, then a decomposition (3.3) of
$\tau_{i+1,j}$ must contain some $t_{i+1,k}$ as well as some $\tau_{h,j}$.
An analogous modification has to be made for the argument following
the decomposition (3.7) of $\tau_{i,j}$.

Let us denote the set of the elements of the root lattice contained in
the cells of $A$ by $\p'$. Since we had to replace $\p$ by $\p'$, the
conclusion of the above arguments is that $t_{i,j}$ is equal to
$$\max \{m:\tau_{i,j}=\b_1+\dots+\b_m,\text{ for some
}\b_l\in\p',\ l=1,2,\dots,m\}.$$
This expression differs from the right-hand side of (3.1) by a
replacement of $\p$ by $\p'$. This indeed makes a difference if the
type that we consider is $B_n$ or $D_n$. 
However, our goal is actually to prove that it
equals exactly the right-hand side of (3.1) as long as $i\ne j+1$, 
for any of the types $B_n$, $C_n$, or $D_n$.  
To see that this is true in type $B_n$, let in a decomposition (3.7)
of $\tau_{i,j}$, where all the $\b_l$'s are in $\p'$, 
one of the $\b_l$'s be equal to a nonroot, to
$\tau_{r,r-1}=2\a_{r}+2\a_{r+1}+\dots+2\a_n$ say. Then some other
$\b_l$ must be necessarily equal to some
$\tau_{s,2n-r+1}=\a_s+\a_{s+1}+\dots+\a_{r-1}$. Then we can replace
these two elements by $\tau_{r,r}$ and $\tau_{s,2n-r}$, both of
which are now roots in $\p$. For, both of which are clearly roots, and
in addition $\tau_{s,2n-r}>\tau_{s,2n-r+1}$, whence $\tau_{s,2n-r}$
must belong to $\p$ since $\tau_{s,2n-r+1}$ does (here we use that
$\p$ is a dual order ideal), and $\tau_{r,r}$ must be in $\p$ because
otherwise we would have $\l_r=r-1$ (that is, the cell $(r,r-1)$ would
be the rightmost cell in row $r$).
An analogous argument holds in type $D_n$.

Now that we have shown that $t_{i,j}$ is equal to the right-hand side
of (3.1), and since also here $t_{1,1}$ is the maximum of all entries
$t_{i,j}$, it follows that $t_{1,1}$ is the class of nilpotence of the
ideal $\i$. Since Proposition~2.1 shows that $t_{1,1}$ can be obtained
by repeated application of (3.8), the proof is complete.
\endemo

In view of this proposition, we have reduced the problem of counting
ideals in $\I$ with a given class of nilpotence to the problem of counting
certain partitions contained in a staircase with respect to the
outcome of the algorithm in Proposition~2.1. 
This partition enumeration problem is already in a convenient form for
being tackled directly in the cases of types $B_n$ and $C_n$. We shall
do that in Sections~6 and 5, respectively. However, in the case of
type $D_n$ we have the additional complication that the partition $\l$
is either given by $\tilde A$ or by $\widetilde {A^\bullet}$. It means that
every partition in which the $(n-1)$st and $n$th column have the same length
is counted once, whereas
every partition whose $(n-1)$st and $n$th column differ in length is
actually counted twice. Therefore we have to distinguish between these
two cases. So, for a Lie algebra $\g$ of type $D_n$,
let, as in Theorem~D.1, 
$\d_n(K)$ denote the number of ideals in $\I$ with class of nilpotence
$K$, let $\d^{(1)}_n(K)$ denote the number of partitions
$\l=(\l_1,\l_2,\dots,\l_{2n-2})$ in the staircase $T_{2n-2}$ with the
property that
$(\l_2,\dots,\l_{2n-2})$ is self-conjugate, that $\l_i\ne i-1$ for $i\ge2$, 
and that $n(\l)=K$, and let
$\d^{(2)}_n(K)$ denote the number of partitions
$\l=(\l_1,\l_2,\dots,\l_{2n-2})$ in the staircase $T_{2n-2}$  
that have the property that
$(\l_2,\dots,\l_{2n-2})$ is self-conjugate, that $\l_i\ne i-1$ for $i\ge2$, 
that the $(n-1)$st and the
$n$th column have the same length, and that $n(\l)=K$. Then the
above arguments show that 
$$\d_n(K)=2\d^{(1)}_n(K)-\d^{(2)}_n(K).\tag4.1$$

In Section~7 we shall show how to compute $\d^{(1)}_n(K)$. 
Very conveniently, the computation of $\d^{(2)}_n(K)$ can be reduced to the
computation of the number of ideals with class of nilpotence $K$ in
type $B_{n-1}$.

\proclaim{Proposition 4.4} For $\g$ a Lie algebra of type $B_n$, 
let, as in Theorem~B.1, $\be_n(K)$ denote the number of 
ideals in $\I$ with class of nilpotence $K$.
Then $\d^{(2)}_n(K)=\b_{n-1}(K)$.\endproclaim
\demo{Proof} By Proposition~4.3 we know that $\be_{n-1}(K)$ is equal
to the number of partitions $\mu=(\mu_1,\mu_2,\dots,\mu_{2n-3})$ 
contained in the
staircase $T_{2n-3}$ with the property that $(\mu_2,\dots,\mu_{2n-3})$ is
self-conjugate, that $\mu_i\ne i-1$ for $i\ge2$, 
and that it satisfies $n(\mu)=K$, where $n(\mu)$ is
the result of the algorithm of Proposition~2.1. 

\midinsert
\vskip10pt
\vbox{
$$
{\Pfad(0,0),222222222111111111\endPfad
\Pfad(0,0),112221221121112212\endPfad
\PfadDicke{.4pt}
\Pfad(0,2),11\endPfad
\Pfad(0,3),11\endPfad
\Pfad(5,6),222\endPfad
\Pfad(6,6),222\endPfad
\SPfad(0,-1),1122211122221111222\endSPfad
\thinlines
\Diagonale(0,-3){12}
\Label\r{x+y=2n-1}(8,3)
\DuennPunkt(2,-1)
\DuennPunkt(5,2)
\DuennPunkt(9,6)
}
\hbox{\hskip7cm}
\Pfad(0,1),2222222211111111\endPfad
\Pfad(0,1),1122122112112212\endPfad
\PfadDicke{.4pt}
\Pfad(0,3),11\endPfad
\Pfad(5,6),222\endPfad
\SPfad(0,0),11222111222111222\endSPfad
\thinlines
\Diagonale(0,-2){11}
\Label\r{x+y=2n-2}(8,4)
\DuennPunkt(2,0)
\DuennPunkt(5,3)
\DuennPunkt(8,6)
\hskip5.5cm
$$
\centerline{\eightpoint a. A partition counted by $\d^{(2)}_n(K)$\hskip1cm
b. A partition counted by $\be_{n-1}(K)$}
\vskip6pt
\centerline{\eightpoint Figure 3}
}
\vskip10pt
\endinsert

Thus it suffices to set up a one-to-one correspondence between the
partitions counted by $\d^{(2)}_n(K)$ and those counted by
$\be_{n-1}(K)$. Such a correspondence is easily set up. Given a
partition counted by $\d^{(2)}_n(K)$ we delete the $(n-1)$st column and
the $n$th row. Clearly, we obtain a partition
$\mu=(\mu_1,\mu_2,\dots,\mu_{2n-3})$ contained in
$T_{2n-3}$ with the property that $(\mu_2,\dots,\mu_{2n-3})$ 
is self-conjugate and that
$\mu_i\ne i-1$ for $i\ge2$. It is equally obvious that the mapping can easily be
reversed. What we need in addition is that the result of 
the algorithm of Proposition~2.1 yields the same value, namely $K$, for the
original partition as well as for the ``reduced" partition. This is
most obvious from the geometric rendering of the algorithm, see
Figure~3, which displays an example with $n=7$ and $K=3$. 
In Part (a) of the figure, the $(n-1)$st column and the $n$th row are
marked, while in Part (b) the places are marked where that column and
that row were cut out. When we delete the
$(n-1)$st column and the $n$th row of the original partition then the
effect on the broken ray is plainly that a horizontal piece of
unit length is
cut out (the one that passed exactly under this column) and that a
vertical piece of unit length is cut out (the one that passed exactly to the
right of that row). Otherwise the ray is completely identical. In
particular, the number of touching points with the bounding line must
be the same, and this number is exactly equal to the outcome of the
algorithm of Proposition~2.1.
\endemo

\heading \S5 Proofs in type $C$\endheading

Let $\g$ be of type $C_n$. We have to prove Theorems~C.2 and C.6, upon
which the other theorems and corollaries labelled C follow, as we
have described in Section~2.

By Proposition~4.3 we know that the number of ideals in $\I$ with class
of nilpotence $K$ is equal to the number of self-conjugate partitions
$\l=(\l_1,\l_2,\dots,\l_{2n-1})$
contained in the staircase $(2n-1,2n-2,\dots,1)$ with $n(\l)=K$, where
$n(\l)$ denotes the result of the algorithm of Proposition~2.1.

The strategy which we use to count the latter partitions
is based on the following observation:
instead of considering the self-conjugate partitions
$(\l_1,\l_2,\dots,\l_{2n-1})$, we consider just
the ``upper halfs," the shifted partitions with row lengths
$(\l_1,\l_2-1,\l_3-2,\dots)$. (In fact, this is the collection of
cells $A$ that contains the roots that define the corresponding ideal;
see Section~4.) Draw the shifted partition
$(\l_1,\l_2-1,\l_2-2,\dots)$, see Figures~4 and 5. Start
as in the $A_{2n-1}$ algorithm. I.e., begin on the vertical edge on
the right of cell $(1,\l_1)$, and move 
downward until we touch the antidiagonal
$x+y=2n$. At the touching point we turn direction from 
vertical-down to horizontal-left, and move on until we touch a vertical
part of the Ferrers diagram. At the touching point we turn direction
from  horizontal-left to vertical-down. We iterate until we meet the
diagonal $x=y$. Let $k$ be the number of touching points on the
antidiagonal $x+y=2n$. Then we claim that 
if we meet the diagonal $x=y$ while travelling
horizontal-left, then the class of nilpotence (of the corresponding
ideal) is $2k$, whereas if we meet the
diagonal while travelling vertical-downward, then the class of
nilpotence is $2k+1$. For example, if $n=9$, then for the shifted partition
$(16,13,11,8,7,5,3)$ in Figure~4 (corresponding to the self-conjugate
partition $(16,14,13,11,11,10,9,6,6,5,4,3,3,2,1,1)$) 
the class of nilpotence is $2\cdot
3=6$, whereas for the partition 
$(16,13,11,8,5,3,1)$ in Figure~5 (corresponding to the self-conjugate
partition $(16,14,13,11,9,8,7,6,5,4,4,3,3,2,1,1)$) 
the class of nilpotence is $2\cdot
2+1=5$. (At this point, the thick segments are without relevance and
should therefore be ignored.)
 
\midinsert
\vskip10pt
\vbox{
$$
\Pfad(-9,9),1111111111111111\endPfad
\Pfad(-3,2),11121212211212112\endPfad
\Pfad(-9,8),2\endPfad
\Pfad(-9,8),1\endPfad
\Pfad(-8,7),2\endPfad
\Pfad(-8,7),1\endPfad
\Pfad(-7,6),2\endPfad
\Pfad(-7,6),1\endPfad
\Pfad(-6,5),2\endPfad
\Pfad(-6,5),1\endPfad
\Pfad(-5,4),2\endPfad
\Pfad(-5,4),1\endPfad
\Pfad(-4,3),2\endPfad
\Pfad(-4,3),1\endPfad
\Pfad(-3,2),2\endPfad
\SPfad(4,7),1112\endSPfad
\SPfad(2,4),1122\endSPfad
\SPfad(-1,1),1122\endSPfad
{\PfadDicke{3pt}
\Pfad(1,3),2\endPfad
\Pfad(4,6),2\endPfad
\Pfad(7,8),2\endPfad}
\thinlines
\Diagonale(-1,-1){11}
\Label\r{x+y=2n}(10,8)
\AntiDiagonale(-10,10){11}
\Label\o{x=y}(-10,10)
\DuennPunkt(1,1)
\DuennPunkt(4,4)
\DuennPunkt(7,7)
\hskip0cm
$$
\centerline{\eightpoint Figure 4}
}
\endinsert

\midinsert
\vbox{
$$
\Pfad(-9,9),1111111111111111\endPfad
\Pfad(-3,2),12121211211212112\endPfad
\Pfad(-9,8),2\endPfad
\Pfad(-9,8),1\endPfad
\Pfad(-8,7),2\endPfad
\Pfad(-8,7),1\endPfad
\Pfad(-7,6),2\endPfad
\Pfad(-7,6),1\endPfad
\Pfad(-6,5),2\endPfad
\Pfad(-6,5),1\endPfad
\Pfad(-5,4),2\endPfad
\Pfad(-5,4),1\endPfad
\Pfad(-4,3),2\endPfad
\Pfad(-4,3),1\endPfad
\Pfad(-3,2),2\endPfad
\SPfad(4,7),1112\endSPfad
\SPfad(-1,4),1111122\endSPfad
\SPfad(-1,1),22\endSPfad
\PfadDicke{3pt}
\Pfad(-1,3),2\endPfad
\Pfad(4,6),2\endPfad
\Pfad(7,8),2\endPfad
\thinlines
\Diagonale(-1,-1){11}
\Label\r{x+y=2n}(10,8)
\AntiDiagonale(-10,10){11}
\Label\o{x=y}(-10,10)
\DuennPunkt(4,4)
\DuennPunkt(7,7)
\hskip0cm
$$
\centerline{\eightpoint Figure 5}
}
\vskip10pt
\endinsert

These claims become immediately obvious from the geometrical
picture. Consider the self-conjugate partition together with the
broken ray that is obtained according to the geometric rendering of
the algorithm of Proposition~2.1. Figure~6.a shows the self-conjugate
partition $(16,14,13,11,11,10,9,6,6,5,4,3,3,2,1,1)$ (corresponding to
the shifted partition in Figure~4) together with the
corresponding broken ray. In the figure
the ray hits the line $x=y$ while travelling
horizontal-left, so that we are in the first of the two possible cases.
If we simply reflect the portion of the ray
which is below the line $x=y$ about this line (see Figure~6.b), 
then it is obvious that
the touching points below the line $x=y$ get reflected just in between
the touching points above $x=y$. Thus, in the first case the total number of
touching points of the ray generated by the original self-conjugate
partition is exactly twice the number of the touching points above
$x=y$, while in the second case it is twice this number plus 1. This
proves our claims.

\def\kleinerKreis(#1,#2){\unskip
  \raise#2 \Einheit\hbox to0pt{\hskip#1 \Einheit
          \raise-4pt\hbox to0pt{\hss$\circ$\hss}\hss}}

\midinsert
\vskip10pt
\vbox{
$$
\Einheit.3cm
\Pfad(-17,-7),22222222222222221111111111111111\endPfad
\Pfad(-17,-7),12212122112121221121212211212112\endPfad
\SPfad(-4,7),1112\endSPfad
\SPfad(-6,4),1122\endSPfad
\SPfad(-9,1),1122\endSPfad
{\PfadDicke{.4pt}
\Pfad(-10,1),1\endPfad
\Pfad(-12,-2),1122\endPfad
\Pfad(-16,-6),1122\endPfad
\Pfad(-17,-8),12\endPfad}
\thinlines
\Diagonale(-18,-10){20}
\Label\r{x+y=2n}(3,8)
\AntiDiagonale(-18,10){11}
\Label\o{x=y}(-18,10)
\kleinerKreis(-16,-8)
\kleinerKreis(-14,-6)
\kleinerKreis(-10,-2)
\DuennPunkt(-7,1)
\DuennPunkt(-4,4)
\DuennPunkt(-1,7)
\hskip0cm
\Pfad(-4,-3),1111111111111111\endPfad
\Pfad(2,-10),11121212211212112\endPfad
\Pfad(-4,-4),2\endPfad
\Pfad(-4,-4),1\endPfad
\Pfad(-3,-5),2\endPfad
\Pfad(-3,-5),1\endPfad
\Pfad(-2,-6),2\endPfad
\Pfad(-2,-6),1\endPfad
\Pfad(-1,-7),2\endPfad
\Pfad(-1,-7),1\endPfad
\Pfad(0,-8),2\endPfad
\Pfad(0,-8),1\endPfad
\Pfad(1,-9),2\endPfad
\Pfad(1,-9),1\endPfad
\Pfad(2,-10),2\endPfad
\SPfad(9,-5),1112\endSPfad
\SPfad(7,-8),1122\endSPfad
\SPfad(4,-11),1122\endSPfad
{\PfadDicke{.4pt}
\Pfad(4,-11),2\endPfad
\Pfad(5,-10),1122\endPfad
\Pfad(9,-6),1122\endPfad
\Pfad(12,-4),12\endPfad}
\thinlines
\Diagonale(4,-13){11}
\Label\r{x+y=2n}(16,-4)
\AntiDiagonale(-5,-2){11}
\Label\o{x=y}(-3,-3)
\DuennPunkt(6,-11)
\DuennPunkt(9,-8)
\DuennPunkt(12,-5)
\kleinerKreis(7,-10)
\kleinerKreis(11,-6)
\kleinerKreis(13,-4)
\hskip0cm
$$
\vskip-1cm
\centerline{\eightpoint a. A self-conjugate partition\hskip6cm}
\vskip1cm
\centerline{\eightpoint \hskip5cm b. One half of the self-conjugate partition}
\vskip6pt
\centerline{\eightpoint Figure 6}
}
\endinsert

We are now in the position to prove Theorem~C.6 (and thereby
Theorem~C.1, which follows by setting $t=1$ in Theorem~C.6). 
Subsequently, we shall prove Theorem~C.2.

\demo{Proof of Theorem C.6} For a
fixed $K$, we aim at
computing the generating function $\sum _{h\ge0} ^{}\gamma_n(h,K)t^h$
for ideals $\i\in\I$ with class of nilpotence $K$. Once this is done,
Equation~(C.6) will follow immediately. Clearly 
the dimension of an ideal $\i$ is equal to the size of the
corresponding shifted partition, i.e., to the sum of its parts. Let us
denote the size of a partition $\mu$ by $\v\mu$. Thus, 
what we are asking for is to compute the generating function $\sum
_{\mu} ^{}t^{\v\mu}$ for shifted partitions $\mu$ contained in the
shifted staircase $S_{2n-1}$ such that the broken ray construction
yields $\fl{K/2}$ touching points on $x+y=2n-1$, and $x=y$ is hit while
travelling horizontal-left if $K$ is even, respectively
vertically-down if $K$ is odd. 

Suppose
we fix a particular broken ray. How do we obtain all the shifted
partitions whose corresponding ray is exactly this fixed one?
The answer is immediate once we observe that every time a broken ray
hits a vertical part of the shifted partition while travelling
horizontal-left, then next the broken ray will continue by a unit step
vertically down. These vertical unit steps are shown as thick segments
in Figures~4 and 5. In turn, every shifted partition which contains
these vertical segments will generate exactly the given broken ray. Let the
$x$-coordinates of these vertical segments be $i_1,i_2,\dots,i_k$, in
ascending order. For example, in Figure~4 we have $i_1=10$, $i_2=13$,
and $i_3=16$, while in Figure~5 we have $i_1=8$, $i_2=13$,
and $i_3=16$.
It is easy to see that $i_1,i_2,\dots,i_k$ determine
a unique broken ray. If we now use the well-known fact that the
generating function $\sum _{\nu} ^{}t^{\v\nu}$,
summed over all (ordinary) partitions $\nu$ which
are contained in an $a\times b$ rectangle, is equal to the
$t$-binomial coefficient
$\left[\smallmatrix {a+b}\\ {b} \endsmallmatrix\right]_t$
(cf\., e.g\., \cite{Sta, Prop.~1.3.19}), then if $K=2k$ we obtain 
$$\multline
\sum _{h\ge0} ^{}\gamma_n(h,2k)t^h=
\sum_{n<i_1<\dots<i_k<i_{k+1}=2n} 
\prod_{j=1}^{k-1}t^{i_{j+1}(i_{j+3}-i_{j+2})}
\bmatrix {i_{j+2}-i_j-1}\\{i_{j+1}-i_j}\endbmatrix_t\\
\cdot
\sum _{\ell=0} ^{i_2-i_1-1}\bmatrix i_1+i_2-2n-1\\\ell\endbmatrix_t
t^{i_1(i_3-i_2)-\binom{2n-i_2+1}2+\binom {\ell+1}2}
\endmultline\tag 5.1$$
(where $i_{k+2}=2n+1$) for the generating function for ideals with
class of nilpotence $2k$, and
$$\multline
\sum _{h\ge0} ^{}\gamma_n(h,2k-1)t^h=
\sum_{2n-i_2<i_1\le n<i_2<\dots<i_k<i_{k+1}=2n} 
t^{i_1(i_3-i_2)-\binom{2n-i_2+1}2}\\
\cdot
(1+t)(1+t^2)\cdots(1+t^{i_1+i_2-2n-1})
\prod_{j=1}^{k-1}t^{i_{j+1}(i_{j+3}-i_{j+2})}
\bmatrix {i_{j+2}-i_j-1}\\{i_{j+1}-i_j}\endbmatrix_t
\endmultline\tag 5.2$$
(where again $i_{k+2}=2n+1$) for the generating function for ideals with
class of nilpotence $2k-1$. 
Because of (C.7), both formulas can be combined.
After having also done the substitution
$i_j\to i_j+n$, the result is (C.6).
\endemo

\demo{Proof of Theorem C.2}
We start by first concentrating on the case that $h$ is even, $h=2k$
say. We want to compute the number of ideals with class of
nilpotence at most $2k$. In principle, we could sum up the expressions
(5.1) and (5.2) with $t=1$, but it is more convenient to use an algorithm
which determines whether the class of nilpotence is {\it at most\/} $2k$ (as
opposed to {\it equal to} $2k$). This
algorithm works as follows. Again we are given a shifted partition
$\l$. This time we start at the diagonal $x=y$ (!),
at the lowest intersection point with the Ferrers diagram of $\l$. We
move right until we touch the antidiagonal $x+y=2n$. At the
touching point we turn direction from horizontal-right to vertical-up,
and move on until we touch a horizontal part of the Ferrers
diagram. At the touching point we turn direction from vertical-up to
horizontal-right. Etc. Then the class of nilpotence is at most twice the
number of touching points on $x+y=2n$. See Figure~7, where 
this procedure is applied to the partition of Figure~4.

\midinsert
\vskip10pt
\vbox{
$$
\Pfad(-9,9),1111111111111111\endPfad
\Pfad(-3,2),11121212211212112\endPfad
\Pfad(-9,8),2\endPfad
\Pfad(-9,8),1\endPfad
\Pfad(-8,7),2\endPfad
\Pfad(-8,7),1\endPfad
\Pfad(-7,6),2\endPfad
\Pfad(-7,6),1\endPfad
\Pfad(-6,5),2\endPfad
\Pfad(-6,5),1\endPfad
\Pfad(-5,4),2\endPfad
\Pfad(-5,4),1\endPfad
\Pfad(-4,3),2\endPfad
\Pfad(-4,3),1\endPfad
\Pfad(-3,2),2\endPfad
\SPfad(7,8),12\endSPfad
\SPfad(3,6),11122\endSPfad
\SPfad(-2,2),1111222\endSPfad
\thinlines
\Diagonale(-1,-1){11}
\Label\r{x+y=2n}(10,8)
\AntiDiagonale(-10,10){11}
\Label\o{x=y}(-10,10)
\DuennPunkt(2,2)
\DuennPunkt(6,6)
\DuennPunkt(8,8)
\hskip0cm
$$
\centerline{\eightpoint Figure 7}
}
\vskip10pt
\endinsert

Let $j_1,j_1+j_2,\dots,j_1+j_2+\dots+j_k$ denote the
$y$-coordinates of the touching points of the broken ray on $x+y=2n$,
in ascending order. For example, in Figure~7 we have $j_1=1$,
$j_2=2$, and $j_3=4$.  
In addition, let us write $j_{k+1}$ for
$n-j_1-j_2-\dots-j_k$. Then,
if we follow the line of argument of the proof of Theorem~C.6 (in a
slightly modified form), we obtain 
$$\underset j_1,j_2,\dots,j_k\ge1,\,j_{k+1}\ge0\to
{\sum _{j_1+j_2+\dots+j_{k+1}=n}} ^{}\binom {j_1+j_2-1}{j_1}
\binom {j_2+j_3-1}{j_2}\cdots \binom {j_{k-1}+j_k-1}{j_{k-1}}
\binom {j_{k}+2j_{k+1}}{j_{k}}$$
for the number of ideals with class of nilpotence $2k$ or $2k-1$, and
therefore, as a moment's thought shows, we obtain
$$\underset j_1,j_2,\dots,j_k,j_{k+1}\ge0\to
{\sum _{j_1+j_2+\dots+j_{k+1}=n}} ^{}\binom {j_1+j_2-1}{j_1}
\binom {j_2+j_3-1}{j_2}\cdots \binom {j_{k-1}+j_k-1}{j_{k-1}}
\binom {j_{k}+2j_{k+1}}{j_{k}}$$
for the number of ideals with class of nilpotence at most $2k$.
Now we can
easily compute the generating function $\sum _{n\ge0}
^{}\gamma^\le_n(2k)x^n$ for the ideals with class of nilpotence at most $2k$. We
substitute the previous expression in the generating function, and
obtain
$$\multline
{\sum _{j_1,j_2,\dots,j_k,j_{k+1}\ge0}} ^{}
x^{j_1+\cdots+j_k+j_{k+1}}
\binom {j_1+j_2-1}{j_1}
\binom {j_2+j_3-1}{j_2}\\
\cdots \binom {j_{k-1}+j_k-1}{j_{k-1}}
\binom {j_{k}+2j_{k+1}}{j_{k}}.
\endmultline$$
Now we perform the summation over $j_1$ and obtain $(1/(1-x))^{j_2}$
in the sum. Next the summation over $j_2$ is performed, etc.,
thus slowly building up continued fractions of the form of the
left-hand side in (5.3) below. If we use the fact that 
$$
\cfrac 1\\1-
\cfrac x\\1-
\cfrac x\\
\cfrac \ddots\\1-x\endcfrac\endcfrac\endcfrac\endcfrac=
\frac{\U_{h}}{\sqrt x\,\U_{h+1}},
\tag 5.3$$
where $x$ occurs $h$ times in the continued fraction 
(again, $\U_k$ is short-hand for $U_k(1/2\sqrt x)$),
then we 
end up with
$$\sum _{j_{k+1}\ge0} ^{}\frac {1} {\sqrt x}\(\frac {\U_k} {\U_{k+1}}\)^{2j_{k+1}+1}.$$
The sum is easily evaluated as it is just a geometric series. This
gives
$$\frac {\U_k\U_{k+1}} 
{\sqrt x\,\(\U_{k+1}^2-\U_{k}^2\)}.$$
Now, it is easy to check that
$$U_k(x)U_{k+1}(x)=U_{2k+1}(x)+U_{2k-1}(x)+\dots+U_1(x)$$
and
$$U_{k+1}^2(x)-U_k^2(x)=U_{2k+2}(x).$$
Hence, our generating function equals
$$\frac {\U_{2k+1}+\U_{2k-1}+\dots+\U_1} 
{\sqrt x\,\U_{2k+2}}.$$
This proves Theorem~C.2 in the case that $h$ is even.
\medskip

Next we compute the generating function for ideals with class of
nilpotence {\it equal to} $2k-1$, i.e., $\sum _{n\ge0}
\gamma_n(2k-1)x^n$, where, as before, $\gamma_n(h)$ is the
number of $ad$-nilpotent ideals in $\I$ with class of nilpotence
equal to $h$. The case $k=1$ has to be done
separately, but this is trivial. (This is already included in the
Kostant--Peterson result.) So let us assume $k\ge2$. If we substitute 
formula (5.2) with $t=1$ in the generating function, then we obtain
$$\multline
\sum _{n\ge0} ^{}
\sum_{2n-i_2<i_1\le n<i_2<\dots<i_k<2n} x^n
\pmatrix {2n-i_{k-1}-1}\\{i_{k}-i_{k-1}}\endpmatrix
\pmatrix {i_{k}-i_{k-2}-1}\\{i_{k-1}-i_{k-2}}\endpmatrix\\
\cdots
\pmatrix {i_{3}-i_1-1}\\{i_{2}-i_1}\endpmatrix
2^{i_1+i_2-2n-1}
\endmultline$$
for our generating function. As before, the sums over
$n,i_k,i_{k-1},\dots,i_3$ are easily computed. On using again that the
resulting continued fractions can be expressed in terms of
Chebyshev polynomials by means of (5.3), and after having replaced $i_1$ by $i_1+n$ and
$i_2$ by $i_2+n$, we obtain 
$$\sum _{-i_2<i_1\le 0<i_2} ^{}x^{i_2+k-1}\(\frac {\U_1}
{x^{k/2-1}\,\U_{k-1}}\)^2 \(\frac {\U_{k-1}} 
{\sqrt x\,\U_k}\)^{i_2-i_1+1}2^{i_1+i_2-1}.$$
Next we sum over $i_1$, and subsequently over $i_2$. In both cases, it
is just geometric series that have to be summed, 
a terminating series when summing over $i_1$
and two nonterminating series when summing over $i_2$. The result,
after cancellation of factors, is
$$\frac {\U_k} {\(\U_k-2\sqrt
x\,\U_{k-1}\)\(\U_k^2-\U_{k-1}^2\)}.$$
It is now routine to verify that this is equal to
$$
\frac {\U_{2k}+\U_{2k-2}+\dots+\U_0} 
{\sqrt x\,\U_{2k+1}}
-\frac {\U_{2k-1}+\U_{2k-3}+\dots+\U_1} 
{\sqrt x\,\U_{2k}}.
$$
We already know that the second term in this difference is the
generating function for ideals with class of nilpotence at most
$2k-2$. Hence the first term must be the generating function for
ideals with class of nilpotence at most $2k-1$.

At this point Theorem~C.2 is completely proved.\endemo

\heading\S6 Proofs in type $B$\endheading

We have to prove Theorem~B.1, upon
which Corollaries B.2 and B.3 follow, as we
have described in Section~2.

\demo{Proof of Theorem~B.1}
By Proposition~4.3 we know that $\b_n(K)$ is equal to the number of
partitions $(\l_1,\l_2,\dots,\l_{2n-1})$ contained in the staircase
$(2n-1,2n-2,\dots,1)$ with the property that $(\l_2,\dots,\l_{2n-1})$
is self-conjugate, that $\l_i\ne i-1$ for $i\ge2$, and that
$n(\l)=K$, where $n(\l)$ stands for
the result of the algorithm of Proposition~2.1. In slight abuse of
terminology, we shall refer to $n(\l)$ as the ``class of nilpotence of
$\l$."

First of all, the cases $K=0$ and $K=1$ can be treated directly,
the case of $K=1$ being contained in the Kostant--Peterson result.

For the proof of the theorem for $K\ge2$, we follow a similar idea as
in the proof of Theorems~C.6 and C.2 in Section~5. Instead of
considering the above (ordinary) partitions
$(\l_1,\l_2,\dots,\l_{2n-1})$, we consider again
just the ``upper halfs," the shifted partitions with row lengths
$(\l_1,\l_2-1,\l_3-2,\dots)$. (In fact, this is the collection of
cells $A$ that contains the roots that define the corresponding ideal;
see Section~4.) Again we need an algorithm which, given the shifted
partition, allows us to find the class of nilpotence $n(\l)$ without
having to go back to the ordinary partition $\l$. The particular
construction that we are going to use requires {\it two} broken
rays (as opposed to just one as in Section~5). These two broken rays
are constructed as follows:
ray~1 starts on the vertical edge on the right of cell
$(1,\l_1)$ (that is, on the vertical right border of
the first row), whereas ray~2 starts on the vertical edge on the
right of cell $(2,\l_2)$ (that is, on the vertical right border of
the second row). Both rays are determined in an algorithmic manner:
starting on that edge, we move downward until we touch the antidiagonal
$x+y=2n$. At the touching point we turn direction from 
vertical-down to horizontal-left, and move on until we touch a vertical
part of the Ferrers diagram. At the touching point we turn direction
from  horizontal-left to vertical-down. Now the procedure is
iterated, until we reach the diagonal $x=y-1$. See
Figures~7--14 for typical examples, with ray~1 (starting at the
right border of the first row) marked as a thin line, and ray~2
(starting at the right border of the second row) marked as a dotted line.

In order to do the computations,
we need to divide the shifted partitions
$(\l_1,\l_2-1,\dots)$ into 7 subclasses. 
To which subclass a particular
partition belongs depends on whether the corresponding 
rays reach $x=y-1$ while
travelling vertical-down or horizontally-left, and whether one ray
is above the other or not, as is detailed below. 
At the same time we shall be able to read off the class of nilpotence
from the broken rays. 

Below we list the 7 subclasses. It is easy to
see that they cover all possibilities. For each subclass we provide
a precise characterization, a typical example, and we describe how to
read off the class of nilpotence.

\medskip
\noindent{\it Case 1}.  A partition belongs to subclass 1, if ray~2 reaches
$x=y-1$ while travelling horizontal-left, and if either ray~1 reaches
$x=y-1$ while travelling horizontal-left weakly below ray~2, after
having touched $x+y=2n$ one more time than ray~2 (see
Figure~8) or ray~1 reaches
$x=y-1$ while travelling vertical-down strictly to the right of ray~2 (see
Figure~9; that some edges are thick is irrelevant at the moment).

\midinsert
\vskip10pt
\vbox{
$$
\Einheit.4cm
\Pfad(-3,2),11212212111221212112\endPfad
\Pfad(-11,11),1111111111111111111\endPfad
\Pfad(-11,10),2\endPfad
\Pfad(-11,10),1\endPfad
\Pfad(-10,9),2\endPfad
\Pfad(-10,9),1\endPfad
\Pfad(-9,8),2\endPfad
\Pfad(-9,8),1\endPfad
\Pfad(-8,7),2\endPfad
\Pfad(-8,7),1\endPfad
\Pfad(-7,6),2\endPfad
\Pfad(-7,6),1\endPfad
\Pfad(-6,5),2\endPfad
\Pfad(-6,5),1\endPfad
\Pfad(-5,4),2\endPfad
\Pfad(-5,4),1\endPfad
\Pfad(-4,3),2\endPfad
\Pfad(-4,3),1\endPfad
\Pfad(-3,2),2\endPfad
\PfadDicke{.4pt}
\Pfad(-1,0),12222111122221111222\endPfad
\SPfad(-2,1),11122222111112222\endSPfad
\PfadDicke{3pt}
\Pfad(-3,2),1\endPfad
\Pfad(0,3),2\endPfad
\Pfad(1,5),2\endPfad
\Pfad(4,7),2\endPfad
\Pfad(6,9),2\endPfad
\Pfad(8,10),2\endPfad
\thinlines
\Diagonale(-1,-1){12}
\Label\r{x+y=2n}(11,9)
\AntiDiagonale(-12,11){12}
\Label\o{x=y-1}(-14,10)
\DuennPunkt(0,0)
\DuennPunkt(1,1)
\DuennPunkt(4,4)
\DuennPunkt(6,6)
\DuennPunkt(8,8)
\Kreis(0,11)
$$
\centerline{\eightpoint Figure 8}
}
\vskip10pt
\endinsert

\midinsert
\vskip10pt
\vbox{
$$
\Einheit.4cm
\Pfad(-3,2),11221212111221212112\endPfad
\Pfad(-11,11),1111111111111111111\endPfad
\Pfad(-11,10),2\endPfad
\Pfad(-11,10),1\endPfad
\Pfad(-10,9),2\endPfad
\Pfad(-10,9),1\endPfad
\Pfad(-9,8),2\endPfad
\Pfad(-9,8),1\endPfad
\Pfad(-8,7),2\endPfad
\Pfad(-8,7),1\endPfad
\Pfad(-7,6),2\endPfad
\Pfad(-7,6),1\endPfad
\Pfad(-6,5),2\endPfad
\Pfad(-6,5),1\endPfad
\Pfad(-5,4),2\endPfad
\Pfad(-5,4),1\endPfad
\Pfad(-4,3),2\endPfad
\Pfad(-4,3),1\endPfad
\Pfad(-3,2),2\endPfad
\PfadDicke{.4pt}
\Pfad(-1,0),22221111122221111222\endPfad
\SPfad(-2,1),11122222111112222\endSPfad
\PfadDicke{3pt}
\Pfad(-3,2),1\endPfad
\Pfad(-1,3),2\endPfad
\Pfad(1,5),2\endPfad
\Pfad(4,7),2\endPfad
\Pfad(6,9),2\endPfad
\Pfad(8,10),2\endPfad
\thinlines
\Diagonale(-1,-1){12}
\Label\r{x+y=2n}(11,9)
\AntiDiagonale(-12,11){12}
\Label\o{x=y-1}(-14,10)
\DuennPunkt(1,1)
\DuennPunkt(4,4)
\DuennPunkt(6,6)
\DuennPunkt(8,8)
\Kreis(0,11)
$$
\centerline{\eightpoint Figure 9}
}
\vskip10pt
\endinsert

Let $k$ be the number of touching points of ray~2 on the antidiagonal
$x+y=2n$. Then
the class of nilpotence of an ideal in this subclass is $2k+1$. 

\smallskip
\noindent{\it Case 2}. A partition belongs to subclass 2, if ray~2 reaches
$x=y-1$ while travelling horizontal-left, and if ray~1 reaches
$x=y-1$ while travelling vertical-down weakly to the left of ray~2 (see
Figure~10).

\midinsert
\vskip10pt
\vbox{
$$
\Einheit.4cm
\Pfad(-3,2),12211212111221212112\endPfad
\Pfad(-11,11),1111111111111111111\endPfad
\Pfad(-11,10),2\endPfad
\Pfad(-11,10),1\endPfad
\Pfad(-10,9),2\endPfad
\Pfad(-10,9),1\endPfad
\Pfad(-9,8),2\endPfad
\Pfad(-9,8),1\endPfad
\Pfad(-8,7),2\endPfad
\Pfad(-8,7),1\endPfad
\Pfad(-7,6),2\endPfad
\Pfad(-7,6),1\endPfad
\Pfad(-6,5),2\endPfad
\Pfad(-6,5),1\endPfad
\Pfad(-5,4),2\endPfad
\Pfad(-5,4),1\endPfad
\Pfad(-4,3),2\endPfad
\Pfad(-4,3),1\endPfad
\Pfad(-3,2),2\endPfad
\PfadDicke{.4pt}
\Pfad(-2,1),22211111122221111222\endPfad
\SPfad(-2,1),11122222111112222\endSPfad
\PfadDicke{3pt}
\Pfad(-3,2),1\endPfad
\Pfad(-2,3),2\endPfad
\Pfad(1,5),2\endPfad
\Pfad(4,7),2\endPfad
\Pfad(6,9),2\endPfad
\Pfad(8,10),2\endPfad
\thinlines
\Diagonale(-1,-1){12}
\Label\r{x+y=2n}(11,9)
\AntiDiagonale(-12,11){12}
\Label\o{x=y-1}(-14,10)
\DuennPunkt(1,1)
\DuennPunkt(4,4)
\DuennPunkt(6,6)
\DuennPunkt(8,8)
\Kreis(0,11)
$$
\centerline{\eightpoint Figure 10}
}
\vskip10pt
\endinsert

Let $k$ be the number of touching points of ray~2 on the antidiagonal
$x+y=2n$. Then
the class of nilpotence of an ideal in this subclass is $2k$.

\smallskip
\noindent{\it Case 3}. A partition belongs to subclass 3, if ray~2 reaches
$x=y-1$ while travelling vertical-down, and if ray~1 reaches
$x=y-1$ while travelling horizontal-left strictly above ray~2 (see
Figure~11).

\midinsert
\vskip10pt
\vbox{
$$
\Einheit.4cm
\Pfad(-4,3),12112211212112112112\endPfad
\Pfad(-11,11),1111111111111111111\endPfad
\Pfad(-11,10),2\endPfad
\Pfad(-11,10),1\endPfad
\Pfad(-10,9),2\endPfad
\Pfad(-10,9),1\endPfad
\Pfad(-9,8),2\endPfad
\Pfad(-9,8),1\endPfad
\Pfad(-8,7),2\endPfad
\Pfad(-8,7),1\endPfad
\Pfad(-7,6),2\endPfad
\Pfad(-7,6),1\endPfad
\Pfad(-6,5),2\endPfad
\Pfad(-6,5),1\endPfad
\Pfad(-5,4),2\endPfad
\Pfad(-5,4),1\endPfad
\Pfad(-4,3),2\endPfad
\PfadDicke{.4pt}
\Pfad(-3,2),11111222222111111222\endPfad
\SPfad(-1,0),22222211111112222\endSPfad
\PfadDicke{3pt}
\Pfad(-4,3),1\endPfad
\Pfad(-1,5),2\endPfad
\Pfad(2,7),2\endPfad
\Pfad(6,9),2\endPfad
\Pfad(8,10),2\endPfad
\thinlines
\Diagonale(-1,-1){12}
\Label\r{x+y=2n}(11,9)
\AntiDiagonale(-12,11){12}
\Label\o{x=y-1}(-14,10)
\DuennPunkt(2,2)
\DuennPunkt(6,6)
\DuennPunkt(8,8)
\Kreis(0,11)
$$
\centerline{\eightpoint Figure 11}
}
\vskip10pt
\endinsert

Let $k$ be the number of touching points of ray~1 on the antidiagonal
$x+y=2n$. Then
the class of nilpotence of an ideal in this subclass is $2k$.

\smallskip
\noindent{\it Case 4}. A partition belongs to subclass 4, if ray~2 reaches
$x=y-1$ while travelling vertical-down, and if ray~1 reaches
$x=y-1$ while travelling vertical-down weakly to the left of ray~2, after
having touched $x+y=2n$ one more time than ray~2 (see
Figure~12).

\midinsert
\vskip10pt
\vbox{
$$
\Einheit.4cm
\Pfad(-4,3),11221211111221212112\endPfad
\Pfad(-11,11),1111111111111111111\endPfad
\Pfad(-11,10),2\endPfad
\Pfad(-11,10),1\endPfad
\Pfad(-10,9),2\endPfad
\Pfad(-10,9),1\endPfad
\Pfad(-9,8),2\endPfad
\Pfad(-9,8),1\endPfad
\Pfad(-8,7),2\endPfad
\Pfad(-8,7),1\endPfad
\Pfad(-7,6),2\endPfad
\Pfad(-7,6),1\endPfad
\Pfad(-6,5),2\endPfad
\Pfad(-6,5),1\endPfad
\Pfad(-5,4),2\endPfad
\Pfad(-5,4),1\endPfad
\Pfad(-4,3),2\endPfad
\PfadDicke{.4pt}
\Pfad(-2,1),22211111122221111222\endPfad
\SPfad(-1,0),22222211111112222\endSPfad
\PfadDicke{3pt}
\Pfad(-4,3),1\endPfad
\Pfad(-2,3),2\endPfad
\Pfad(-1,5),2\endPfad
\Pfad(4,7),2\endPfad
\Pfad(6,9),2\endPfad
\Pfad(8,10),2\endPfad
\thinlines
\Diagonale(-1,-1){12}
\Label\r{x+y=2n}(11,9)
\AntiDiagonale(-12,11){12}
\Label\o{x=y-1}(-14,10)
\DuennPunkt(4,4)
\DuennPunkt(6,6)
\DuennPunkt(8,8)
\Kreis(0,11)
$$
\centerline{\eightpoint Figure 12}
}
\vskip10pt
\endinsert

Let $k$ be the number of touching points of ray~1 on the antidiagonal
$x+y=2n$. Then
the class of nilpotence of an ideal in this subclass is $2k$.

\smallskip
\noindent{\it Case 5}. A partition belongs to subclass 5, if ray~2 reaches
$x=y-1$ while travelling vertical-down, and if ray~1 reaches
$x=y-1$ while travelling horizontal-left weakly below ray~2 (see
Figure~13).

\midinsert
\vskip10pt
\vbox{
$$
\Einheit.4cm
\Pfad(-4,3),12221112121112112112\endPfad
\Pfad(-11,11),1111111111111111111\endPfad
\Pfad(-11,10),2\endPfad
\Pfad(-11,10),1\endPfad
\Pfad(-10,9),2\endPfad
\Pfad(-10,9),1\endPfad
\Pfad(-9,8),2\endPfad
\Pfad(-9,8),1\endPfad
\Pfad(-8,7),2\endPfad
\Pfad(-8,7),1\endPfad
\Pfad(-7,6),2\endPfad
\Pfad(-7,6),1\endPfad
\Pfad(-6,5),2\endPfad
\Pfad(-6,5),1\endPfad
\Pfad(-5,4),2\endPfad
\Pfad(-5,4),1\endPfad
\Pfad(-4,3),2\endPfad
\PfadDicke{.4pt}
\Pfad(-2,1),11122222221111111222\endPfad
\SPfad(-3,2),22221111111112222\endSPfad
\PfadDicke{3pt}
\Pfad(-4,3),1\endPfad
\Pfad(-3,5),2\endPfad
\Pfad(1,7),2\endPfad
\Pfad(6,9),2\endPfad
\Pfad(8,10),2\endPfad
\thinlines
\Diagonale(-1,-1){12}
\Label\r{x+y=2n}(11,9)
\AntiDiagonale(-12,11){12}
\Label\o{x=y-1}(-14,10)
\DuennPunkt(1,1)
\DuennPunkt(6,6)
\DuennPunkt(8,8)
\Kreis(0,11)
$$
\centerline{\eightpoint Figure 13}
}
\vskip10pt
\endinsert

Let $k$ be the number of touching points of ray~1 on the antidiagonal
$x+y=2n$. Then
the class of nilpotence of an ideal in this subclass is $2k-1$.

\smallskip
\noindent{\it Case 6}. A partition belongs to subclass 6, if ray~2 reaches
$x=y-1$ while travelling vertical-down, and if ray~1 reaches
$x=y-1$ while travelling vertical-down weakly to the right of ray~2,
both rays touching $x+y=2n$ an equal number of times (see
Figure~14).

\midinsert
\vskip10pt
\vbox{
$$
\Einheit.4cm
\Pfad(-4,3),12222121111112112112\endPfad
\Pfad(-11,11),1111111111111111111\endPfad
\Pfad(-11,10),2\endPfad
\Pfad(-11,10),1\endPfad
\Pfad(-10,9),2\endPfad
\Pfad(-10,9),1\endPfad
\Pfad(-9,8),2\endPfad
\Pfad(-9,8),1\endPfad
\Pfad(-8,7),2\endPfad
\Pfad(-8,7),1\endPfad
\Pfad(-7,6),2\endPfad
\Pfad(-7,6),1\endPfad
\Pfad(-6,5),2\endPfad
\Pfad(-6,5),1\endPfad
\Pfad(-5,4),2\endPfad
\Pfad(-5,4),1\endPfad
\Pfad(-4,3),2\endPfad
\PfadDicke{.4pt}
\Pfad(-2,1),22222221111111111222\endPfad
\SPfad(-3,2),22221111111112222\endSPfad
\PfadDicke{3pt}
\Pfad(-4,3),1\endPfad
\Pfad(-3,5),2\endPfad
\Pfad(-2,7),2\endPfad
\Pfad(6,9),2\endPfad
\Pfad(8,10),2\endPfad
\thinlines
\Diagonale(-1,-1){12}
\Label\r{x+y=2n}(11,9)
\AntiDiagonale(-12,11){12}
\Label\o{x=y-1}(-14,10)
\DuennPunkt(6,6)
\DuennPunkt(8,8)
\Kreis(0,11)
$$
\centerline{\eightpoint Figure 14}
}
\vskip10pt
\endinsert

Let $k-1$ be the number of touching points of ray~1 on the antidiagonal
$x+y=2n$. Then
the class of nilpotence of an ideal in this subclass is $2k-1$.

\smallskip
\noindent{\it Case 7}. A partition belongs to subclass 7, if ray~2 reaches
$x=y-1$ while travelling horizontal-left, and if ray~1 reaches
$x=y-1$ while travelling horizontal-left weakly above ray~2, both rays
touching $x+y=2n$ an equal number of times (see
Figure~15). 

\midinsert
\vskip10pt
\vbox{
$$
\Einheit.4cm
\Pfad(-5,4),11121112111221212112\endPfad
\Pfad(-11,11),1111111111111111111\endPfad
\Pfad(-11,10),2\endPfad
\Pfad(-11,10),1\endPfad
\Pfad(-10,9),2\endPfad
\Pfad(-10,9),1\endPfad
\Pfad(-9,8),2\endPfad
\Pfad(-9,8),1\endPfad
\Pfad(-8,7),2\endPfad
\Pfad(-8,7),1\endPfad
\Pfad(-7,6),2\endPfad
\Pfad(-7,6),1\endPfad
\Pfad(-6,5),2\endPfad
\Pfad(-6,5),1\endPfad
\Pfad(-5,4),2\endPfad
\PfadDicke{.4pt}
\Pfad(-5,4),11111111122221111222\endPfad
\SPfad(-2,1),11122222111112222\endSPfad
\PfadDicke{3pt}
\Pfad(-5,4),1\endPfad
\Pfad(1,5),2\endPfad
\Pfad(4,7),2\endPfad
\Pfad(6,9),2\endPfad
\Pfad(8,10),2\endPfad
\thinlines
\Diagonale(-1,-1){12}
\Label\r{x+y=2n}(11,9)
\AntiDiagonale(-12,11){12}
\Label\o{x=y-1}(-14,10)
\DuennPunkt(1,1)
\DuennPunkt(4,4)
\DuennPunkt(6,6)
\DuennPunkt(8,8)
\Kreis(0,11)
$$
\centerline{\eightpoint Figure 15}
}
\vskip10pt
\endinsert

Let $k$ be the number of touching points of ray~2 on the antidiagonal
$x+y=2n$. Then
the class of nilpotence of an ideal in this subclass is $2k$.

\medskip
Let $\b^{(i)}_n(K)$ denote the number of partitions in the $i$th
subclass with class of nilpotence $K$, $i=1,2,\dots,7$. We will now,
for each subclass, compute the corresponding generating function $\sum
_{n\ge0} ^{}\b^{(i)}_n(K)x^n$. This is made possible by formulas for
$\b^{(i)}_n(K)$ for each subclass in the spirit of (A.1).

\medskip
\noindent
{\it The enumeration in Case 1}. 
For given $n$, the number $\b^{(1)}_n(2k+1)$ of partitions in this
subclass is given by the multiple summation
$$\sum _{-i_2\le i_1\le i_2\le \dots\le i_{2k+1}\le n-1} ^{}
\bigg(\prod _{j=4} ^{2k+2}\binom {i_j-i_{j-1}+i_{j-2}-i_{j-3}-1}
{i_{j-2}-i_{j-3}}\bigg)\kern-4pt \sum _{\ell=0} ^{i_3-i_2-1}\binom {i_3+i_1-1}\ell,$$
where, as always in the sequel, $i_{2k+2}=n-1$. This expression is
obtained from the geometrical presentation as in Figures~7 or 8, by
denoting the deviation of the first vertical edges on each downward
travel along a ray (in Figures~7 and 8 these are the thick vertical
edges; the thick horizontal edge marks the downmost and leftmost edge
of the shifted Ferrers diagram; it must touch the diagonal $x=y-1$) from the
``reference point" $(n,0)$ (in Figures~7 and 8 it is marked by a
circle) by $i_{2k+1},i_{2k},\dots,i_1$, from right to left. Thus, in
Figure~8 we have $i_5=8$, $i_4=6$, $i_3=4$, $i_2=1$, $i_1=0$, and in
Figure~9 we have $i_5=8$, $i_4=6$, $i_3=4$, $i_2=1$, $i_1=-1$. Then
arguments very similar to those that we used in the proofs of
Theorems~C.6 and C.2 in the previous section show
that, for fixed $i_{2k+1},i_{2k},\dots,i_1$, the number of
partitions is equal to the summand in the above sum.

We now compute the generating function $\sum _{n\ge0}
^{}\be_n^{(1)}(2k+1)x^n$. By definition, this is
$$\multline \sum _{n\ge0} ^{}x^n
\sum _{-i_2\le i_1\le i_2\le \dots\le i_{2k+1}\le n-1} ^{}
\bigg(\prod _{j=4} ^{2k+2}\binom {i_j-i_{j-1}+i_{j-2}-i_{j-3}-1}
{i_{j-2}-i_{j-3}}\bigg)\\
\cdot \sum _{\ell=0} ^{i_3-i_2-1}\binom {i_3+i_1-1}\ell.
\endmultline$$
Now we would like to perform the sum over $n$, then the sum over
$i_{2k+1}$, then over $i_{2k}$, etc., one after the other. In order to
do this conveniently, we have to split the range of the sum into the
possibilities 
$$-i_2\le i_1\le i_2,\quad i_2+\ell+1\le i_3< \dots< i_{2k+1}< n-1$$
and
$$\multline
0\le i_1=i_2,\quad i_2+\ell+1\le i_3=i_4<i_5=i_6<
\cdots\\
<i_{2j-1}=i_{2j}<i_{2j+1}< \dots< i_{2k+1}< n-1,
\endmultline$$
for some $j$ between $2$ and $k+1$. These two do indeed cover all
possibilities which contribute to the sum. For, if we have
$i_s=i_{s-1}$ for some $s\ge4$, then, in order that the summand does
not vanish, we must also have $i_{s-2}=i_{s-3}$, etc. The case of $s$
being even is covered by the second possibility above. If $s$ should
be odd, then we would also have $i_3=i_2$. Since then the sum over
$\ell$ is empty, this does not contribute to the sum.

Now, in each of the two cases the sums over $n$, $i_{2k+1}$, \dots,
$i_4$ are easily carried out by the binomial theorem. Thereby,
continued fractions as in the left-hand side of (5.3)
are slowly built up.

To finish the calculation, we use the identity
$$\sum _{\ell=0} ^{i_3-i_2-1}\binom {i_3+i_1-1}\ell=
\sum _{\ell=0} ^{i_3-i_2-1}\binom {i_1+i_2+\ell}\ell\cdot
\cases 2^{i_3-i_2-\ell-2}&\text{if }\ell<i_3-i_2-1,\\
1&\text{if }\ell=i_3-i_2-1.\endcases\tag6.1$$
For, now the sums over $i_3$, $\ell$, $i_1$, and finally $i_2$, in
that order, can easily be carried out, as this amounts again to just
summing binomial or even just geometric series. To get rid of the
continued fractions, we use Equation~(5.3).
The result is
$${\frac{{\sqrt{x}}\,\U_{k+1}\,
      \left( {{\U_{k}}^2} + {{\U_{k+1}}^2} \right) }{{{\U_{k}}^2}\,
   {{\U_{2k+2}}^2}   \left(   \U_{k}-2\,{\sqrt{x}}\,\U_{k-1} \right) \,
      }}.
$$

The computation for the second possible range is similar. The result
is
$${\frac{{\sqrt{x}}\,{{\U_{k+1}}^2}}
    { \U_{k-j+1}\,\U_{k-j+2}\,
\U_{k}\, \U_{2k+2}\left(   \U_{k}-2\,{\sqrt{x}}\,\U_{k-1} \right) \,
    }}.
$$

In order to obtain the overall generating function for the ideals in
this subclass, we have to take the first expression and add to it the
sum of the second expressions over $j$ from $2$ to $k+1$. Using the
identity
$$\frac {1} {U_0(x)U_1(x)}+\frac {1} {U_1(x)U_2(x)}+\dots+\frac {1}
{U_s(x)U_{s+1}(x)}=\frac {U_s(x)} {U_{s+1}(x)},$$
we obtain
$$ {\frac{{\sqrt{x}}\,\U_{k+1}\,
      \left( {{\U_{k}}^2} + {{\U_{k+1}}^2} + 
        \U_{k-1}\,\U_{k+1}\,\U_{2k+2} \right) }{{{\U_{k}}^2}\,
   {{\U_{2k+2}}^2}   \left(   \U_{k}-2\,{\sqrt{x}}\,\U_{k-1} \right) \,
      }}.
$$

\smallskip
\noindent
{\it The enumeration in Case 2}. By arguments analogous to the ones
in Case~1, for given $n$ the number $\b^{(2)}_n(2k)$ of partitions in this
subclass is given by the multiple summation
$$\underset -i_3+1\le i_1\le -i_2-1\to
{\sum _{0\le i_2\le \dots\le i_{2k+1}\le n-1} ^{}}
\bigg(\prod _{j=4} ^{2k+2}\binom {i_j-i_{j-1}+i_{j-2}-i_{j-3}-1}
{i_{j-2}-i_{j-3}}\bigg) 2^{i_3+i_1-1}.$$

We now compute the generating function $\sum _{n\ge0}
^{}\be_n^{(2)}(2k)x^n$. 
We proceed as in Case 1. Here, the analysis is in fact easier. First,
there is no sum over $\ell$, so that we do not have to use the
transformation formula (6.1). Second, both $i_2=i_1$ and $i_3=i_2$
produce vanishing summands, so that we can restrict the summation to
$$-i_3+1\le i_1\le -i_2-1,\quad 0\le i_2< \dots< i_{2k+1}< n-1.$$
As a result, we obtain for the generating function the expression
$${\frac{{\sqrt{x}}\,\U_{k-1}\,\U_{k+1}}
    {{{\U_{k}}^2}\,\U_{2k+2} \,
      \left(   \U_{k+1}-\U_{k-1} \right)\,
   \left(   \U_{k}-2\,{\sqrt{x}}\,\U_{k-1} \right) }}.
$$

\smallskip
\noindent
{\it The enumeration in Case 3}. By arguments analogous to the ones
in Case~1, for given $n$ the number $\b^{(3)}_n(2k)$ of partitions in this
subclass is given by the multiple summation
$$\sum _{-i_2\le i_1<0\le i_2\le \dots\le i_{2k}\le n-1} ^{}
\kern-2pt
\bigg(\prod _{j=4} ^{2k+1}\binom {i_j-i_{j-1}+i_{j-2}-i_{j-3}-1}
{i_{j-2}-i_{j-3}}\bigg)\kern-4pt \sum _{\ell=0} ^{i_3-i_2-1}\binom
{i_3+i_1-1}\ell .$$

We now compute the generating function $\sum _{n\ge0}
^{}\be_n^{(3)}(2k)x^n$. 
We proceed as in Case 1. Since both $i_2=i_1$ and $i_3=i_2$
produce vanishing summands, we can restrict the summation to
$$-i_2\le i_1<0\le i_2< \dots< i_{2k}< n-1.$$
As a result, we obtain for the generating function the expression
$${\frac{{\sqrt{x}}\,\U_{k-1}\,\U_{k+1}}
    {{{\U_{k}}^2}\, \U_{2k}\,\left(   \U_{k+1}-\U_{k-1} \right)
\left(   \U_{k}-2\,{\sqrt{x}}\,\U_{k-1} \right) \,
      }}.
$$

\smallskip
\noindent
{\it The enumeration in Case 4}. By arguments analogous to the ones
in Case~1, for given $n$ the number $\b^{(4)}_n(2k)$ of partitions in this
subclass is given by the multiple summation
$$
{\sum _{-i_3+1\le i_1\le i_2<0\le i_3\le \dots\le i_{2k+1}\le n-1} ^{}}
\bigg(\prod _{j=4} ^{2k+2}\binom {i_j-i_{j-1}+i_{j-2}-i_{j-3}-1}
{i_{j-2}-i_{j-3}}\bigg) 2^{i_3+i_1-1}.$$

We now compute the generating function $\sum _{n\ge0}
^{}\be_n^{(4)}(2k)x^n$. 
We proceed as in Case 1. We split the range of summation into the
possibilities 
$$-i_3+1\le i_1\le i_2<0\le i_3< \dots< i_{2k+1}< n-1$$
and
$$ -i_3+1\le i_1= i_2<0\le i_3=i_4< \cdots
<i_{2j-1}=i_{2j}<i_{2j+1}<\dots<i_{2k+1}< n-1,
$$
for some $j$ between $2$ and $k+1$.
The contribution to the generating function of the first range is
$${\frac{{\sqrt{x}}\,{{\U_{k-1}}^2}}
    {{{\U_{k}}^2}\, \U_{2k}\,\left(   \U_{k+1}-\U_{k-1} \right) 
\left(   \U_{k}-2\,{\sqrt{x}}\,\U_{k-1} \right) \,
     }},
$$
while the contribution of the second range is
$${\frac{{\sqrt{x}}\,{{\U_{k-1}}^2}}
    {\U_{k-j+1}\,\U_{k-j+2}\,\U_{k}\,\U_{2k}
\left(   \U_{k}-2\,{\sqrt{x}}\,\U_{k-1} \right) \,
      }}.
$$
Summing the second expression over $j$ from $2$ to $k+1$, and adding
the result to the first expression gives
$${\frac{{\sqrt{x}}\,{{\U_{k-1}}^2}\,
      }
      {{{\U_{k}}^2}\, \U_{2k}
\left(   \U_{k}-2\,{\sqrt{x}}\,\U_{k-1} \right) \,
     }}
 \left(   {\frac{1}{  \U_{k+1}-\U_{k-1}}}+\U_{k-1} \right).
$$

\smallskip
\noindent
{\it The enumeration in Case 5}. By arguments analogous to the ones
in Case~1, for given $n$ the number $\b^{(5)}_n(2k-1)$ of partitions in this
subclass is given by the multiple summation
$$\underset -i_3+1\le i_1\le -i_2-1\to
{\sum _{0\le i_2\le \dots\le i_{2k}\le n-1} ^{}}
\bigg(\prod _{j=4} ^{2k+1}\binom {i_j-i_{j-1}+i_{j-2}-i_{j-3}-1}
{i_{j-2}-i_{j-3}}\bigg) 2^{i_3+i_1-1}.$$

We now compute the generating function $\sum _{n\ge0}
^{}\be_n^{(5)}(2k-1)x^n$. 
We proceed as in Case 1. Since both $i_2=i_1$ and $i_3=i_2$
produce vanishing summands, we can restrict the summation to
$$ -i_3+1\le i_1\le -i_2-1,\quad 0\le i_2< \dots< i_{2k}< n-1$$
As a result, we obtain for the generating function the expression
$${\frac{{\sqrt{x}}\,\U_{k-1}}
    {{{\U_{2k}}^2}\left(   \U_{k}-2\,{\sqrt{x}}\,\U_{k-1} \right) }
    }.
$$

\smallskip
\noindent
{\it The enumeration in Case 6}. By arguments analogous to the ones
in Case~1, for given $n$ the number $\b^{(6)}_n(2k-1)$ of partitions in this
subclass is given by the multiple summation
$$
{\sum _{-i_3+1\le i_1\le i_2<0\le i_3\le \dots\le i_{2k}\le n-1} ^{}}
\bigg(\prod _{j=4} ^{2k+1}\binom {i_j-i_{j-1}+i_{j-2}-i_{j-3}-1}
{i_{j-2}-i_{j-3}}\bigg) 2^{i_3+i_1-1}.$$

We now compute the generating function $\sum _{n\ge0}
^{}\be_n^{(6)}(2k-1)x^n$. 
We proceed as in Case 1. We split the range of summation into the
possibilities 
$$-i_3+1\le i_1\le i_2<0\le i_3< \dots< i_{2k}< n-1$$
and
$$\multline 
-i_3+1\le i_1= i_2<0\le i_3=i_4<\cdots\\
<i_{2j-1}=i_{2j}<i_{2j+1}< \dots< i_{2k}< n-1,
\endmultline$$
for some $j$ between $2$ and $k$.
The contribution to the generating function of the first range is
$${\frac{{\sqrt{x}}\,\U_{k-1}}
    {{{\U_{2k}}^2}\left(   \U_{k}-2\,{\sqrt{x}}\,\U_{k-1} \right) }
    },
$$
while the contribution of the second range is
$${\frac{{\sqrt{x}}\,{{\U_{k-1}}^2}}
    {\U_{k-j}\,\U_{k-j+1}\,\U_{k}\, \U_{2k}
\left(   \U_{k}-2\,{\sqrt{x}}\,\U_{k-1} \right) \,
     }}.
$$
Summing the second expression over $j$ from $2$ to $k$, and adding
the result to the first expression gives
$${\frac{{\sqrt{x}}\,\U_{k-1}\,
     }{{{\U_{2k}}^2}
      \left(   \U_{k}-2\,{\sqrt{x}}\,\U_{k-1} \right) 
      }}
 \left( 1 + {\frac{\U_{k-2}\,\U_{2k}}{\U_{k}}} \right) .
$$

\smallskip
\noindent
{\it The enumeration in Case 7}. By arguments analogous to the ones
in Case~1, for given $n$ the number $\b^{(7)}_n(2k)$ of partitions in this
subclass is given by the multiple summation
$$\sum _{0\le i_1\le i_2\le \dots\le i_{2k}\le n-1} ^{}
\bigg(\prod _{j=4} ^{2k+1}\binom {i_j-i_{j-1}+i_{j-2}-i_{j-3}-1}
{i_{j-2}-i_{j-3}}\bigg) \sum _{\ell=0} ^{i_3-i_2-1}\binom {i_3+i_1-1}\ell.$$

We now compute the generating function $\sum _{n\ge0}
^{}\be_n^{(7)}(2k)x^n$. 
We proceed as in Case 1. We split the range of summation into the
possibilities 
$$0\le i_1\le i_2< \dots< i_{2k}< n-1$$
and
$$0\le i_1= i_2<i_3=i_4<\dots<i_{2j-1}<i_{2j}<i_{2j+1}< \dots< i_{2k}<
n-1,$$
for some $j$ between $2$ and $k$.
The contribution to the generating function of the first range is
$${\frac{{\sqrt{x}}\,{{\U_{k+1}}^3}}
    {\U_{k-1}\,{{\U_{k}}^2}\, \U_{2k+2}
 \left(   \U_{k+1}-\U_{k-1} \right)
      \left(   \U_{k}-2\,{\sqrt{x}}\,\U_{k-1} \right) \,
     }},$$
while the contribution of the second range is
$$ {\frac{{\sqrt{x}}\,{{\U_{k+1}}^2}}
    { \U_{k-j}\,\U_{k-j+1}\,\U_{k}\,\U_{2k+2}
\left(   \U_{k}-2\,{\sqrt{x}}\,\U_{k-1} \right) \,
     }}.
$$
Summing the second expression over $j$ from $2$ to $k$, and adding
the result to the first expression gives
$$ {\frac{{\sqrt{x}}\,{{\U_{k+1}}^2}\,
     }{
      \U_{k-1}\,{{\U_{k}}^2}\,\U_{2k+2}
      \left(   \U_{k}-2\,{\sqrt{x}}\,\U_{k-1} \right) }} \left( 
        {\frac{\U_{k+1}}{  \U_{k+1}-\U_{k-1}}}+
 \U_{k-2}\,\U_{k}  \right) .
$$

\medskip
Finally, in order to complete the proof of Theorem~B.1, in case that
$K$ is even, the expressions of Cases~2, 3, 4 and 7 have to be summed,
while in case that $K$ is odd, the expressions of Cases~1 (with $k$ replaced
by $k-1$), 5 and 6 have to be summed. Some simplification yields the
claimed expressions.
\endemo

\heading\S7 Proofs in type D\endheading

We have to prove Theorem~D.1, upon
which Corollaries D.2 and D.3 follow, as we
have described in Section~2. 

\demo{Proof of Theorem~D.1}
According to (4.1) and Proposition~4.4 we have to compute 
$$2\sum _{n\ge0} ^{}\d^{(1)}_n(K)x^n-\sum _{n\ge0}
^{}\b_{n-1}(K)x^n.\tag7.1$$
By Theorem~B.1, we already know an expression for the second sum. It
remains to compute $\sum _{n\ge0} ^{}\d^{(1)}_n(K)x^n$, where
$\d^{(1)}_n(K)$ is the number of partitions
$\l=(\l_1,\l_2,\dots,\l_{2n-2})$ in the staircase $T_{2n-2}$ with the
property that
$(\l_2,\dots,\l_{2n-2})$ is self-conjugate, that $\l_i\ne i-1$ for $i\ge2$, 
and that $n(\l)=K$.

We proceed by imitating the arguments
used in Case $B$ in the previous section. First, instead of
considering the above (ordinary) partitions
$(\l_1,\l_2,\dots,\l_{2n-2})$, we consider again
just the ``upper halfs," the shifted partitions with row lengths
$(\l_1,\l_2-1,\l_3-2,\dots)$.
We apply again the construction of the two broken
rays as in Case~B, with the slight modification that it is now the
line $x+y=2n-1$ (instead of $x+y=2n$) where the rays get reflected.
We divide the shifted partitions that we consider here 
again into 7 subclasses. The characterization of the 7 subclasses and
the description of how to determine the corresponding class of
nilpotence are identical to those in Case~B, except that, of course, 
again the line $x+y=2n$ has to be replaced by the line $x+y=2n-1$. 
We therefore omit to repeat them here and instead refer the reader to
Section~6. For writing down formulas for the number of partitions in
each subclass we also follow the derivations in Case~B. In
particular, we choose again $(n,0)$ as the reference point with
respect to which deviations of the first vertical edges on each downward
travel along a ray are measured. 

We now discuss each of the 7 subclasses.
Since everything is completely parallel to the computations in
Case~B, given in the previous section, we can be brief here. For each
subclass, we provide the formula in the spirit of (A.1)
for the number of shifted partitions in that subclass, and the
corresponding generating function.

\medskip
{\it The enumeration in Case 1}: 
For given $n$, the number of shifted partitions in this subclass 
is given by the multiple summation
$$\sum _{-i_2-1\le i_1\le i_2\le \dots\le i_{2k+1}\le n-2} ^{}
\bigg(\prod _{j=4} ^{2k+2}\binom {i_j-i_{j-1}+i_{j-2}-i_{j-3}-1}
{i_{j-2}-i_{j-3}}\bigg)\kern-4pt \sum _{\ell=0} ^{i_3-i_2-1}\binom
{i_3+i_1}\ell,$$
where $i_{2k+2}=n-1$.
For the range 
$$-i_2-1\le i_1\le i_2,\quad i_2+\ell+1\le i_3< \dots< i_{2k+1}< n-2$$
we obtain
$$\frac {2\,x\,\U_{k+1}^2} {\U_{k} \U_{2k+2}^2
   \left(\U_{k}-2\sqrt{x}\U_{k-1}\right)},
$$
for the corresponding generating function,
while for the range 
$$\multline
0\le i_1=i_2,\quad i_2+\ell+1\le i_3=i_4<i_5=i_6<
\cdots\\
<i_{2j-1}=i_{2j}<i_{2j+1}< \dots< i_{2k+1}< n-2,
\endmultline$$
for some $j$ between $2$ and $k+1$, we obtain
$$\frac {x\, \U_{k+1}} {\U_{k-j+1}\,\U_{k-j+2}\,\U_{2k+2}\,
   \left(\U_{k}-2\sqrt{x}\U_{k-1}\right)}.
$$
The overall generating function for this subclass is
$$\frac {x\, \U_{k+1}\left(2\U_{k+1}+\U_{k-1}\U_{2k+2}\right)} 
 {\U_{k}\,\U_{2k+2}^2\,
   \left(\U_{k}-2\sqrt{x}\U_{k-1}\right)
  }.
$$
\medskip\noindent
{\it The enumeration in Case 2}: 
For given $n$, the number of shifted partitions in this subclass 
is given by the multiple summation
$$\underset -i_3\le i_1\le -i_2-2\to
{\sum _{0\le i_2\le \dots\le i_{2k+1}\le n-2} ^{}}
\bigg(\prod _{j=4} ^{2k+2}\binom {i_j-i_{j-1}+i_{j-2}-i_{j-3}-1}
{i_{j-2}-i_{j-3}}\bigg) 2^{i_3+i_1},$$
where $i_{2k+2}=n-1$.
We can restrict the sum to the range
$$-i_3\le i_1\le -i_2-2,\quad 0\le i_2< \dots< i_{2k+1}< n-2.$$
We obtain
$$\frac {x\, \U_{k-1}}
{\U_{k}\,\U_{2k+2}\,\left(\U_{k+1}-\U_{k-1}\right)\,
   \left(\U_{k}-2\sqrt{x}\U_{k-1}\right)}
$$
for the corresponding generating function.

\medskip\noindent
{\it The enumeration in Case 3}: 
For given $n$, the number of shifted partitions in this subclass 
is given by the multiple summation
$$\sum _{-i_2-1\le i_1<0\le i_2\le \dots\le i_{2k}\le n-2} ^{}
\bigg(\prod _{j=4} ^{2k+1}\binom {i_j-i_{j-1}+i_{j-2}-i_{j-3}-1}
{i_{j-2}-i_{j-3}}\bigg)\kern-4pt \sum _{\ell=0} ^{i_3-i_2-1}\binom
{i_3+i_1}\ell,$$
where $i_{2k+1}=n-1$.
We can restrict the sum to the range
$$-i_2-1\le i_1<0\le i_2< \dots< i_{2k}< n-2.$$
We obtain
$$\frac {x\, \U_{k+1}}
{\U_{k}\,\U_{2k}\,\left(\U_{k+1}-\U_{k-1}\right)\,
   \left(\U_{k}-2\sqrt{x}\U_{k-1}\right)}.
$$
for the corresponding generating function.

\medskip\noindent
{\it The enumeration in Case 4}: 
For given $n$, the number of shifted partitions in this subclass 
is given by the multiple summation
$$
\sum _{-i_3\le i_1\le i_2<0\le i_3\le \dots\le i_{2k+1}\le n-2}
\bigg(\prod _{j=4} ^{2k+2}\binom {i_j-i_{j-1}+i_{j-2}-i_{j-3}-1}
{i_{j-2}-i_{j-3}}\bigg) 2^{i_3+i_1},$$
where $i_{2k+2}=n-1$.
For the range 
$$-i_3\le i_1\le i_2<0\le i_3< \dots< i_{2k+1}< n-2$$
we obtain
$$\frac {x\, \U_{k-1}}
{\U_{k}\,\U_{2k}\,\left(\U_{k+1}-\U_{k-1}\right)\,
   \left(\U_{k}-2\sqrt{x}\U_{k-1}\right)},
$$
for the corresponding generating function,
while for the range
$$ -i_3\le i_1= i_2<0\le i_3=i_4< \cdots
<i_{2j-1}=i_{2j}<i_{2j+1}<\dots<i_{2k+1}< n-2,
$$
for some $j$ between $2$ and $k+1$,
we obtain
$$\frac {x\, \U_{k-1}} {\U_{k-j+1}\,\U_{k-j+2}\,\U_{2k}\,
   \left(\U_{k}-2\sqrt{x}\U_{k-1}\right)}.
$$
The overall generating function is
$$\frac {x\, \U_{k-1}} 
  {\U_{k}\,\U_{2k}\,
   \left(\U_{k}-2\sqrt{x}\U_{k-1}\right)}
\left(\frac {1} {\U_{k+1}-\U_{k-1}}+\U_{k-1}\right).
$$

\medskip\noindent
{\it The enumeration in Case 5}: 
For given $n$, the number of shifted partitions in this subclass 
is given by the multiple summation
$$\underset -i_3\le i_1\le -i_2-2\to
{\sum _{0\le i_2\le \dots\le i_{2k}\le n-2} ^{}}
\bigg(\prod _{j=4} ^{2k+1}\binom {i_j-i_{j-1}+i_{j-2}-i_{j-3}-1}
{i_{j-2}-i_{j-3}}\bigg) 2^{i_3+i_1},$$
where $i_{2k+1}=n-1$.
We can restrict the sum to the range
$$ -i_3\le i_1\le -i_2-2,\quad 0\le i_2< \dots< i_{2k}< n-2.$$
We obtain
$$\frac {x\, \U_{k-1}^2} {\U_{2k}^2\,\U_{k}\,
   \left(\U_{k}-2\sqrt{x}\U_{k-1}\right)}
$$
for the corresponding generating function.

\medskip\noindent
{\it The enumeration in Case 6}: 
For given $n$, the number of shifted partitions in this subclass 
is given by the multiple summation
$$
\sum _{-i_3\le i_1\le i_2<0\le i_3\le \dots\le i_{2k}\le n-2}
\bigg(\prod _{j=4} ^{2k+1}\binom {i_j-i_{j-1}+i_{j-2}-i_{j-3}-1}
{i_{j-2}-i_{j-3}}\bigg) 2^{i_3+i_1},$$
where $i_{2k+1}=n-1$.
For the range 
$$-i_3\le i_1\le i_2<0\le i_3< \dots< i_{2k}< n-2$$
we obtain
$$\frac {x\, \U_{k}} {\U_{2k}^2\,
   \left(\U_{k}-2\sqrt{x}\U_{k-1}\right)},
$$
for the corresponding generating function,
while for the range
$$\multline 
-i_3+1\le i_1= i_2<0\le i_3=i_4<\cdots\\
<i_{2j-1}=i_{2j}<i_{2j+1}< \dots< i_{2k}< n-1,
\endmultline$$
for some $j$ between $2$ and $k$, we obtain
$$\frac {x\, \U_{k-1}} {\U_{k-j+1}\,\U_{k-j}\,\U_{2k}\,
   \left(\U_{k}-2\sqrt{x}\U_{k-1}\right)}.
$$
The overall generating function is
$$\frac {x\, \U_{k}} {\U_{2k}^2\,
   \left(\U_{k}-2\sqrt{x}\U_{k-1}\right)}
\left(1+\frac {\U_{k-2}\U_{2k}} {\U_{k}}\right).
$$

\medskip\noindent
{\it The enumeration in Case 7}: 
For given $n$, the number of shifted partitions in this subclass 
is given by the multiple summation
$$\sum _{0\le i_1\le i_2\le \dots\le i_{2k+1}\le n-2} ^{}
\bigg(\prod _{j=4} ^{2k+2}\binom {i_j-i_{j-1}+i_{j-2}-i_{j-3}-1}
{i_{j-2}-i_{j-3}}\bigg)\kern-4pt \sum _{\ell=0} ^{i_3-i_2-1}\binom
{i_3+i_1}\ell,$$
where $i_{2k+2}=n-1$.
For the range
$$0\le i_1\le i_2< \dots< i_{2k}< n-2$$
we obtain
$$\frac {x\, \U_{k+1}^2}
{\U_{k}\,\U_{k-1}\,\U_{2k+2}\,\left(\U_{k+1}-\U_{k-1}\right)\,
   \left(\U_{k}-2\sqrt{x}\U_{k-1}\right)},
$$
for the corresponding generating function,
while for the range
$$0\le i_1= i_2<i_3=i_4<\dots<i_{2j-1}<i_{2j}<i_{2j+1}< \dots< i_{2k}<
n-2,$$
for some $j$ between $2$ and $k$, we obtain
$$\frac {x \,\U_{k+1}} {\U_{k-j+1}\,\U_{k-j}\,\U_{2k+2}\,
   \left(\U_{k}-2\sqrt{x}\U_{k-1}\right)}.
$$
The overall generating function is
$$\frac {x\,\U_{k+1}} 
{\U_{k-1}\,\U_{k}\,\U_{2k+2}\,
   \left(\U_{k}-2\sqrt{x}\U_{k-1}\right)}
\left(\frac {\U_{k+1}} {\U_{k+1}-\U_{k-1}}+\U_{k-2}\U_{k}\right).
$$

\medskip
Finally, in order to complete the proof of Theorem~D.1, we have to
combine all our results to obtain the generating function (7.1).
In order to obtain the generating
function $\sum _{n\ge0} ^{}\d^{(1)}_n(K)x^n$, in case that
$K$ is even, the expressions of Cases~2, 3, 4 and 7 have to be summed,
while in case that $K$ is odd, the expressions of Cases~1 (with $k$ replaced
by $k-1$), 5 and 6 have to be summed. The result is then substituted
in (7.1), together with the result of the previous section
for $\sum _{n\ge0}^{}\b_{n-1}(K)x^n$. Some simplification eventually
yields the claimed expressions.
\endemo

\heading\S8 The exceptional types\endheading

It is not difficult to write down a computer program which determines explicitly the descending central series 
of a given ideal. The final results are given in Table~1.

\midinsert
$$\smatrix
\format\qquad\c\quad&\qquad\c\qquad&\qquad\c\qquad&\qquad\c\qquad&
\qquad\c\qquad &\qquad\c\qquad\\
n(\i)&E_6&E_7&E_8&F_4 &G_2\\
\hlinefor6\\
\hphantom{0}0&\hphantom{00}1&\hphantom{000}1&\hphantom{0000}1&\hphantom{0}1&1\\
\hphantom{0}1&\hphantom{0}63&\hphantom{0}127&\hphantom{00}255 &15&3\\
\hphantom{0}2&210&\hphantom{0}662&\hphantom{0}2200 &28&2\\
\hphantom{0}3&217&\hphantom{0}894&\hphantom{0}3804 &21&1\\
\hphantom{0}4&150&\hphantom{0}766&\hphantom{0}3872 &14&0\\
\hphantom{0}5&\hphantom{0}92&\hphantom{0}576&\hphantom{0}3372 &12&1\\
\hphantom{0}6&\hphantom{0}51&\hphantom{0}403&\hphantom{0}2752 &\hphantom{0}5\\
\hphantom{0}7&\hphantom{0}28&\hphantom{0}279&\hphantom{0}2182 &\hphantom{0}4\\
\hphantom{0}8&\hphantom{0}12&\hphantom{0}175&\hphantom{0}1656 &\hphantom{0}2\\
\hphantom{0}9&\hphantom{00}6&\hphantom{0}115&\hphantom{0}1277 &\hphantom{0}2\\
10&\hphantom{00}2&\hphantom{00}68&\hphantom{00}955 &\hphantom{0}0\\
11&\hphantom{00}1&\hphantom{00}44&\hphantom{00}737&\hphantom{0}1\\
12& &\hphantom{00}23&\hphantom{00} 536&\\
13& &\hphantom{00}14&\hphantom{00}412 &\\
14& &\hphantom{000}7&\hphantom{00}300 &\\
15& &\hphantom{000}4&\hphantom{00}227 &\\
16& &\hphantom{000}1&\hphantom{00}157 &\\
17& &\hphantom{000}1&\hphantom{00}123 &\\
18& & &\hphantom{000}81&\\
19& & &\hphantom{000}61 &\\
20& & &\hphantom{000}40&\\
21& & &\hphantom{000}30 &\\
22& & &\hphantom{000}18 &\\
23& & &\hphantom{000}14 &\\
24& & &\hphantom{0000}7 &\\
25& & &\hphantom{0000}5&\\
26& & &\hphantom{0000}3 &\\
27& & &\hphantom{0000}2 &\\
28& & &\hphantom{0000} 0&\\
29& & &\hphantom{0000}1&\\
\hlinefor6\\
\text{total}\\
\text{number}&833& 4160& 25080& 105&8
\endsmatrix$$
\vskip6pt
\centerline{\eightpoint Table 1}
\endinsert

Just a few remarks on how to get the data in the table. 
Clearly we can work at the level of the root system. If $\i\in \I$ is
encoded  by $\p\subseteq\Dp$, then $n(\i)$ equals the maximal integer $k$ such that    $\theta\in \p^k$ (here
$\theta$ is the highest root of $\D$ and $\p^k$ is inductively defined as $\p^1=\p,\ \p^k=(\p^{k-1}+\p)\cap \D$).
The input for the program (i.e., the $\p$'s) can be  reduced to the determination of the  antichains 
of the root poset. (Indeed, for any finite poset, there is 
a canonical bijection mapping the antichain $\{a_1,\ldots,a_k\}$ to
the dual  order ideal which is the union of the principal dual order ideals
$V_{a_1},\ldots,V_{a_k}$.)  
In turn, the calculation of the antichains has  been done using the
Maple program {\it Coxeter} \cite{Ste}.
Note that the total number of ideals in type $E_6$ is 833,  one more than
the number given in \cite{Sh, Theorem~3.6}. 
Our counting agrees with formula (1.1).

\enddocument